\documentclass{amsproc}
\pdfoutput=1

\newtheorem{theorem}{Theorem}[section]
\newtheorem{lemma}[theorem]{Lemma}
\newtheorem{proposition}[theorem]{Proposition}

\theoremstyle{definition}
\newtheorem{definition}[theorem]{Definition}
\newtheorem{example}[theorem]{Example}

\theoremstyle{remark}
\newtheorem{remark}[theorem]{Remark}
\newtheorem{notation}[theorem]{Notation}
\newtheorem{constr}[theorem]{Construction}

\numberwithin{equation}{section}

\newcommand{\dash}{{\operatorname{-}}}

\usepackage{enumerate}
\usepackage{mathtools}
\usepackage{tikz-cd}
\usepackage[all]{xy}
\usepackage{graphicx, adjustbox}
\usepackage{extpfeil}
\usepackage{comment}
\usepackage{float}
\usepackage[labelformat=empty]{caption}
\usepackage[toc]{appendix}
\usepackage{hyperref}
\usepackage{cleveref}
\usepackage{resizegather}
\usepackage[activate={true, nocompatibility}, final, tracking=true, kerning=true, spacing=true, factor=1100, stretch=10, shrink=10]{microtype}
\usepackage{xcolor}
\usepackage[normalem]{ulem}
\usepackage{pbox}

\DeclareMathOperator{\Bun}{Bun}
\DeclareMathOperator{\Aff}{Aff}
\DeclareMathOperator{\Spec}{Spec}
\DeclareMathOperator{\Sym}{Sym}
\DeclareMathOperator{\Iso}{Iso}
\DeclareMathOperator{\Red}{Red}
\DeclareMathOperator{\rk}{rk}

\DeclareMathOperator{\Gr}{Gr}
\DeclareMathOperator{\Hom}{Hom}
\DeclareMathOperator{\Coh}{Coh}

\DeclareMathOperator{\Pic}{Pic}

\DeclareMathOperator{\GL}{GL}

\DeclareMathOperator{\QCoh}{QCoh}

\newcommand{\cF}{{\mathcal{F}}}
\newcommand{\cO}{{\mathcal{O}}}
\newcommand{\cM}{{\mathcal{M}}}
\newcommand{\inj}{\hookrightarrow}

\makeatletter
\newcommand{\colim@}[2]{%
  \vtop{\m@th\ialign{##\cr
    \hfil$#1\operator@font colim$\hfil\cr
    \noalign{\nointerlineskip\kern1.5\ex@}#2\cr
    \noalign{\nointerlineskip\kern-\ex@}\cr}}%
}
\newcommand{\colim}{%
  \mathop{\mathpalette\colim@{\rightarrowfill@\textstyle}}\nmlimits@
}
\makeatother

\makeatletter
\DeclareRobustCommand
  \myvdots{\vbox{\baselineskip4\p@ \lineskiplimit\z@
    \hbox{.}\hbox{.}\hbox{.}}}
\makeatother

\tikzset{
  symbol/.style={
    draw=none,
    every to/.append style={
      edge node={node [sloped, allow upside down, auto=false]{$#1$}}}
  }
}
\makeatletter
\tikzset{
  column sep/.code=\def\pgfmatrixcolumnsep{\pgf@matrix@xscale*(#1)},
  row sep/.code   =\def\pgfmatrixrowsep{\pgf@matrix@yscale*(#1)},
  matrix xscale/.code=%
    \pgfmathsetmacro\pgf@matrix@xscale{\pgf@matrix@xscale*(#1)},
  matrix yscale/.code=%
    \pgfmathsetmacro\pgf@matrix@yscale{\pgf@matrix@yscale*(#1)},
  matrix scale/.style={/tikz/matrix xscale={#1},/tikz/matrix yscale={#1}}}
\def\pgf@matrix@xscale{1}
\def\pgf@matrix@yscale{1}
\makeatother

\begin{document}

\title{A guide to moduli theory beyond GIT}


\author[T. L. G\'omez]{Tom\'as L. G\'omez}

\address{Instituto de Ciencias Matem\'aticas (CSIC-UAM-UC3M-UCM),
Nicol\'as Cabrera 15, Campus Cantoblanco UAM, 28049 Madrid, Spain}

\email{tomas.gomez@icmat.es}

\author[A. Fernández Herrero]{Andres Fernández Herrero}

\address{Mathematics Hall MC 4404, 2990 Broadway,
New York, N.Y. 10027, USA}

\email{af3358@columbia.edu}

\author[A. Zamora]{Alfonso Zamora}

\address{Departamento de Matem\'atica Aplicada a las TIC, ETSI Inform\'aticos, Universidad Polit\'ecnica de Madrid, Campus de Montegancedo,
28660 Madrid, Spain}

\email{alfonso.zamora@upm.es}

\thanks{This work is supported
by grants CEX2019-000904-S and PID2019-108936GB-C21 (funded by MCIN/AEI/ 10.13039/501100011033).}

%    General info
\date{February 1 2023}

\dedicatory{Dedicated to Peter E. Newstead, on his $80^{th}$ birthday}


\begin{abstract}
In this survey we provide an overview of some recent developments in the construction of moduli spaces using stack-theoretic techniques. We will also explain the analogue of Harder-Narasimhan stratifications for general stacks, known as $\Theta$-stratifications. As an application of the ideas exposed here, we address the moduli problem of principal bundles over higher dimensional projective varieties, as well as its different compactifications by the so-called principal $\rho$-sheaves. We construct a stratification by instability types whose lower strata admits a proper good moduli space of ``Gieseker semistable" objects and a new Gieseker-type Harder-Narasimhan filtration for these objects. Detailed proofs of the latter results will appear elsewhere.
\end{abstract}

\maketitle

\section{Introduction}

The construction and study of moduli spaces has been one of the driving forces in algebraic geometry for decades. After the algebraic construction of the Picard scheme by Grothendieck, the study of the moduli of vector bundles came next. Not only has Peter Newstead been one of the pioneers in this study since his early papers in the 60's, but by founding the group VBAC  (Vector Bundles on Algebraic Curves), he has been pivotal in building a large and successful community dedicated to this and related areas. It is an honour for us to dedicate this survey to him.

Since its inception by Mumford, Geometric Invariant Theory (GIT) \cite{mumford-git} has been a powerful way to construct moduli spaces for many moduli problems. When performing GIT, we do not attempt to classify all possible objects in a given moduli problem, but rather only those declared to be semistable. After adding the additional data of a ``framing", these objects are parametrized by a scheme called the parameter space. Mumford's GIT provides a way to define the quotient of this scheme by the action of an algebraic group encoding the framing, revisiting Hilbert's $14^{th}$ problem, to obtain a moduli scheme parametrizing semistable objects. 

The intrinsic point of view in the Beyond GIT program developed in \cite{halpernleistner2021structure} does not focus on the subclass of semistable objects but on the whole moduli stack. A generalization of the Hilbert-Mumford criterion leads to a new stacky numerical criterion through test stacks taking the role of former one-parameter subgroups. This makes the notion of semistability naturally appear in the stack, avoiding the tedious checking that semistable objects correspond to GIT semistable ones and vice-versa. 
Thanks to the recent paper \cite{alper2019existence}, there are some intrinsic conditions that can be checked locally in an algebraic stack in order to show that it admits a separated good moduli space in the sense of Alper \cite{alper-good-moduli}.

Unstable objects usually come in a natural way with a maximal degeneration (a stacky Harder-Narasimhan filtration). The theory of $\Theta$-stratifications in \cite{halpernleistner2021structure} can be used to find such maximal degenerations without having to consider non-canonically a GIT problem sufficiently large to contain the required unstable object as in \cite{gomez-sols-zamora}. This opens the door to understanding Harder-Narasimhan filtrations also in an intrinsic way, directly from the whole stack. 

Harder-Narasimhan filtrations (also known as canonical filtrations) were first considered in \cite{harder-narasimhan} for vector bundles on curves. Shatz studied its behaviour in families in \cite{shatz}. In higher dimensions we have Gieseker semistablity, and we have Gieseker-Harder-Narasimahn filtrations for torsion free sheaves. Nitsure used them to produce schematic stratifications in \cite{nitsurehnsheaves} (see also \cite{gurjar-nitsure} for $\Lambda$-modules), which give rise to algebraic
stacks corresponding to each Harder-Narasimhan type, thus defining a filtration of
the stack of pure coherent sheaves by locally closed substacks.

As an application of this theory, we address the moduli problem of principal bundles over higher dimensional projective varieties, as well as its different compactifications by the so-called principal $\rho$-sheaves. We construct in \cite{rho-sheaves} a stratification (the $\Theta$-stratification of a given numerical invariant), whose lower open strata admit proper good moduli spaces (as in \cite{alper-good-moduli}) of ``Gieseker semistable'' objects \cite{gieseker_torsion_free}. This gives an intrinsic stack-theoretic construction to the already known moduli space of principal bundles on higher dimensional varieties \cite{schmitt.singular, gomezsols.principalsheaves, glss.singular.char, glss.large}, as well as its compactifications.

For principal bundles, a Harder-Narasimhan filtration is usually called canonical reduction. They were announced by Ramanathan in \cite{ramanathan-copenhagen} and have been constructed over curves by Behrend \cite{behrend-thesis} and Biswas-Holla \cite{biswas-holla-hnreduction}. A reduction to a parabolic subgroup is also defined in \cite{atiyah-bott} by a different method, and it is proved in  \cite{anchouche-hassan-biswas} to coincide with the canonical filtration. Gurjar and Nitsure shows that this reduction also gives a schematic stratification (i.e., a stratification of the stack) \cite{gurjar-nitsure,nitsuregurjar2}. 
For principal bundles in positive characteristics, Gurjar and Nitsure
\cite{gurjar2020hardernarasimhan}
show that there are algebraic stacks corresponding to each HN-type,
and these are radicial
over the algebraic stack of all principal bundles. If the Behrend conjecture
holds then these stacks are shown to define a stratification of the
stack of all principal bundles (see also \cite{heinloth-behrends-conjecture}).
However, over higher dimension varieties, where a polynomial Gieseker stability condition is used to construct the compactification of the moduli space, canonical reductions for principal bundles in the literature only take into account the leading coefficients of the reduced Hilbert polynomial (analogous to slope stability as described in \cite{mumford-projetive-invariants}). There has been partial attempts to construct a filtration that accounts for all of the terms in the Hilbert polynomial in lower rank \cite{zamora-rank2-tensors}. On the other hand, in \cite{biswas-zamora-ghn-principal} it was shown that the notion of Gieseker Harder-Narasimhan parabolic reduction cannot be related to the Gieseker-Harder-Narasimhan filtration of the underlying vector bundles for faithful representations, which meant that the main strategies in the literature (e.g. \cite{anchouche-hassan-biswas}) could not possibly apply if we try to account for all terms of the Hilbert polynomial.
 
An important outcome of the results in \cite{rho-sheaves} is a new notion of a Gieseker-Harder-Narasimhan filtration for principal $\rho$-sheaves, which is constructed by iterating a filtration coming from the $\Theta$-stratification (called the leading term HN filtration). The Gieseker Harder-Narasimhan filtration thus obtained consists of a parabolic reduction of the underlying rational principal bundle that measures instability by considering all degree terms of the numerical polynomials. 
Detailed proofs of the results contained in this survey will appear elsewhere.

\begin{notation}
For this survey we work over a fixed ground field $k$, which we assume to be of characteristic $0$ unless explicitly stated. We denote by $\text{Sch}_k$ the category of schemes over $k$. We will sometimes work instead over the smaller category $\text{Aff}_k$ of affine schemes over $k$. For any $k$-schemes $X$ and $S$, we will write $X_S$ to denote the fiber product $X\times S$, which is a product over $k$ unless otherwise stated. 
\end{notation}

\begin{section}{Intrinsic approach to moduli problems}
\begin{subsection}{Algebraic stacks}
Moduli theory rests on the insight that many classification problems in algebraic geometry can be parametrized by schemes. An example of this is the projective space $\mathbb{P}^n_k$, which classifies line bundles along with the choice of $n+1$ generating sections. In order to naturally deal with the algebraic structure of such classification problems, the Grothendieck school of algebraic geometry initiated the functorial point of view on such moduli problems.

 We can view any scheme $X$ as a functor $h_X$ from $(\text{Aff}_k)^{op}$ into sets via the Yoneda embedding. We say that functors of the form $h_X$ are representable by schemes. We will equip $\text{Sch}_k$ and $\text{Aff}_k$ with the additional structure of the Grothendieck topology generated by fppf (finitely presented faithfully flat) covers. All representable functors satisfy descent with respect to this fppf topology.

Equipped with the functorial point of view, Grothendieck proved the representability of some of the central moduli problems in algebraic geometry: the Quot functor and the Picard functor. We refer the reader to the articles by Nitsure and Kleiman in \cite{fga-explained} for a nice exposition of these achievements. While the definition of the Quot functor intuitively captures the sought after moduli problem, surprises arise when one tries to parametrize line bundles on a given projective variety $X$. The natural functor classifying isomorphism classes of line bundles does not satisfy fppf descent, and therefore cannot be represented by a scheme. However, Grothendieck proved that the fppf sheafification $\text{Pic}(X)$ of the natural functor is indeed represented by a scheme.

When the variety $X$ is proper but not projective, Artin showed the fppf sheafification $\text{Pic}(X)$ is represented by an algebraic space \cite[Thm 7.3]{artin-algebrization-formal-moduli-i}. Algebraic spaces form a subcategory of all fppf sheaves. This subcategory provides a more robust setting to parametrize moduli problems, while remaining fundamentally geometric at the same time. We refer to \cite{knutson-algebraic-spaces}, \cite{stacks-project} for a detailed treatments of the theory of algebraic spaces. Artin's result on the representability of $\Pic(X)$ is one of the many achievements that hinted at the usefulness of enlarging the category of schemes in order to study moduli problems.

A more serious issue arises when one tries to generalize the Picard scheme to the setting of higher rank vector bundles on a given projective variety $X$. Even in the case when $X$ is a curve, it becomes impossible to parametrize the moduli problem of all vector bundles in terms of a separated scheme. The main feature of this moduli problem is the existence of large nontrivial automorphisms of the objects, which makes the condition of descent for isomorphism classes a more delicate affair. One approach, pioneered by Mumford when $X$ is a smooth projective curve, is to restrict to a special class of reasonably ``rigid'' vector bundles called stable vector bundles \cite{mumford-projetive-invariants}. Alternatively, if one wants to be able to parametrize all vector bundles, then it becomes necessary to account for their automorphisms in order to naturally deal with the descent condition. This is elegantly encoded in the notion of an algebraic stack, which is another geometric enlargement of the category of schemes originally introduced by Deligne and Mumford \cite{deligne-mumford-irreducibility}.

For the rest of the survey we will assume that the reader has some familiarity with the theory of algebraic stacks. For an informal introduction we refer the reader to \cite{gomez-algebraic-stacks} \cite{heinloth-stack-vector-bundles}; on the other hand we recommend \cite{casalaina-martin-wise-intro-stacks} \cite{alper-notes} for more detailed developments of the theory. Some standard references are \cite{lmb-champsalgebriques}, \cite{olsson-algebraic-stacks} and \cite{stacks-project}.

We will favor the point of view of stacks as pseudofunctors valued in groupoids, which we believe is more intuitive for moduli problems. There should be no issue in translating things into the theoretically cleaner language of fibered categories thanks to \cite[Part 1, \S 3.1.3]{fga-explained}.

Recall that a stack $\mathcal{M}$ is a pseudofuctor from $(\text{Sch}_k)^{op}$ into groupoids that satisfies fppf descent. We say that a stack is algebraic if it admits a smooth morphism $f: U \to \mathcal{M}$ representable by algebraic spaces such that the source $U$ is a scheme. Such morphism $f$ is called an atlas. Many of the properties of algebraic stacks (locally of finite type, quasi-compact, smooth, normal) can be encoded in terms of properties of some atlas. 

Just as in the case of schemes, algebraic stacks are closed under the operation of fiber product. One can define the diagonal of a stack in a similar way as in the case of schemes. The main difference is that in the case of algebraic stacks the diagonal is not a monomorphism. Instead, the fibers of the diagonal encode automorphisms of the corresponding points of the stack.

Throughout this survey we will only work with stacks that are locally of finite type over $k$. Moreover we shall restrict our attention to stacks with affine diagonal. 
\end{subsection}
\begin{subsection}{Good moduli spaces}
Many moduli problems in algebraic geometry, such as the moduli of vector bundles on a smooth projective curve, are parametrized by algebraic stacks \cite{heinloth-stack-vector-bundles}. As in the case of the Picard scheme or the substack of stable vector bundles, it is often possible to ``rigidify" the functor to remove the automorphisms, thus obtaining a functor valued in sets that is representable by an algebraic space. Such algebraic space is called the moduli space associated to the stack.

When we are working with a stack $\mathcal{M}$ with finite inertia, the theorem of Keel and Mori \cite{keel-mori-groupoids} \cite[\href{https://stacks.math.columbia.edu/tag/0DUT}{Tag 0DUT}]{stacks-project} provides a general procedure to obtain a ``rigidification" in the form of an algebraic space $M$ called the coarse moduli space as in the following.

\begin{definition}
    A morphism $f:\cM \to M$ into an algebraic space $M$ is called a coarse moduli space morphism (and $M$ is called a coarse moduli space) if the following are satisfied
    \begin{enumerate}[(1)]
        \item $f$ is initial among morphisms to algebraic spaces.
        \item For any algebraically closed field $K$, the natural transformation $f:\cM(K) \to M(K)$ induces a bijection on isomorphism classes of points. 
    \end{enumerate}
\end{definition}

There is a coarse moduli space morphism $\mathcal{M} \to M$ that relates both functors. Many geometric properties of the original algebraic stack $\mathcal{M}$ (and of the moduli problem parametrized by it) are reflected in the algebraic space $M$, which lives in a closer setting to the more geometric world of varieties.

When the objects parametrized by the algebraic stack $\mathcal{M}$ have infinite automorphisms (i.e. algebraic automorphism groups of positive dimension), then the power of the Keel-Mori theorem becomes unavailable. Unfortunately, many natural moduli problems arising in algebraic geometry have infinite symmetry groups, so it becomes important to develop technical tools to deal with this type of situation. Mumford's Geometric Invariant Theory (GIT) \cite{mumford-git} offers a powerful device to construct a moduli space in the case when the algebraic stack $\mathcal{M}$ is a quotient stack of the form $\left[X/G\right]$, where $X$ is open inside a projective-over-affine variety and $G$ is a reductive group. In this case, the properties of the corresponding morphism $\left[X/G\right] \to M$ are captured by the notion of a ``good quotient".

Alper's notion of good moduli space is an elegant way to encode and generalize Mumford's theory to the setting of arbitrary algebraic stacks that do not necesarily arise as quotient stacks.
\begin{definition}[Good moduli space, {\cite[Defn. 4.1]{alper-good-moduli}}] \label{defn: good moduli space} A quasi-compact and quasiseparated morphism $\varphi: \mathcal{M} \to M$ from an algebraic stack $\mathcal{M}$ to an algebraic space $M$ is called a good moduli space if both of the following properties are satisfied:
\begin{enumerate}[(a)]
\item The pushforward morphism $\varphi_*: \QCoh(\mathcal{M}) \to \text{QCoh}(M)$ is exact.
\item The induced morphism of structure sheaves $\mathcal{O}_{M} \to \varphi_*\mathcal{O}_{\mathcal{M}}$ is an isomorphism.
\end{enumerate}
We call the algebraic space $M$ the good moduli space of $\mathcal{M}$.
\end{definition}

\begin{remark}
The notion of good moduli space becomes quite restrictive if we work over a base of positive or mixed characteristic. Most of the classical moduli spaces that admit a good moduli space in characteristic $0$ (such as the moduli of semistable vector bundles on a curve) do not admit a good moduli space over a field of positive characteristic. Instead, they admit morphisms $\mathcal{M} \to M$ to an algebraic space satisfying a weakening of condition (a) in Definition \ref{defn: good moduli space}. A more appropriate concept, the notion of adequate moduli space, was introduced by Alper in \cite{alper-adequate} to approach moduli problems in arbitrary characteristic.
\end{remark}

The above Definition \ref{defn: good moduli space} is deceivingly concise; it has a lot of strong implications which link the properties of the moduli space $M$ and the stack $\mathcal{M}$. The first main property justifies the name ``the" moduli space.
\begin{theorem}[{\cite[Thm. 6.6]{alper-good-moduli}}]
Let $\mathcal{M}$ be a locally Noetherian algebraic stack with a good moduli space $\phi: \mathcal{M} \to M$. Then $\phi$ is initial among all maps from $\mathcal{M}$ into algebraic spaces. In particular the good moduli space $\phi: \mathcal{M} \to M$ is unique up to unique isomorphism.
\end{theorem}
Another important property of good moduli spaces is that they behave well with respect to base-change.
\begin{proposition}[{\cite[Prop. 4.7(i)]{alper-good-moduli}}]
Let $\mathcal{M}$ be an algebraic stack with a good moduli space $\phi: \mathcal{M} \to M$. For any morphism of algebraic spaces $X \to M$, the base-change $\phi_X: \mathcal{M}\times_{M} X \to X$ is a good moduli space for the fiber product $\mathcal{M} \times_{M} X$.
\end{proposition}

We finally summarize some of the important properties of good moduli spaces in the following
\begin{proposition}
Let $\phi: \mathcal{M} \to M$ be a good moduli space. Then the following hold.
\begin{enumerate}[(i)]
\item $\phi$ is surjective and universally closed \cite[Thm. 4.16(i)+(ii)]{alper-good-moduli}.
\item If $\phi$ has affine diagonal, then $\phi$ is of finite type \cite[Thm. A.1]{alper-hall-rydh-etale-slice}.
\item The fibers of $\phi$ are geometrically connected, each $\phi$-fiber contains a unique closed point \cite[Thm. 4.16(vii) + Prop. 9.1]{alper-good-moduli}.
\item If $\mathcal{M}$ is locally Noetherian, then so is $M$, and the push-forward $\phi_*$ preserves coherence of sheaves \cite[Thm. 4.16(x)]{alper-good-moduli}.
\item If $\mathcal{M}$ is of finite type over an excellent scheme $S$, then $M$ is also of finite type over $S$ \cite[Thm. 4.16(xi)]{alper-good-moduli}.
\end{enumerate}
\end{proposition}

In addition, there are some results that relate the singularities of a stack and the singularities of its good moduli space. We would like to highlight the following recent generalization of a classical theorem of Boutot.
\begin{theorem}[{\cite[Thm. 5]{boutot-thm-good-moduli}}]
Let $\mathcal{M}$ be a smooth algebraic stack with affine diagonal over a field $k$ with good moduli space $\phi: \mathcal{M} \to M$. If $M$ is quasiprojective, then $M$ is of klt type.
\end{theorem}

We end this section by providing some examples of good moduli spaces.
\begin{example}[Affine modulo linearly reductive] \label{example: good moduli space for affine over reductive}
Let $X = \text{Spec}(A)$ be an affine variety equipped with the action of a linearly reductive group $G$ over $k$. We consider the quotient stack $\left[X/G\right]$. We denote by $A^G \subset A$ the subring of invariants. There is a morphism $\phi: \left[X/G\right] \to \text{Spec}(A^G)$ which exhibits $\text{Spec}(A^G)$ as a good moduli space for $\left[X/G\right]$. Indeed, let us check both properties of a good moduli space morphism.
\begin{enumerate}[(a)]
\item The abelian category of quasicoherent sheaves $\text{QCoh}(\left[X/G\right])$ on the stack $\left[X/G\right]$ is equivalent to the category of $G$-equivariant $A$-modules. The push-forward $\phi_*(-)$ takes a $G$-equivariant $A$-module to the module of $\phi_*(M) = M^G \subset M$ of $G$-invariants, viewed as a module over $A^G$. Since the algebraic group $G$ is assumed to be linearly reductive, the functor of taking $G$-invariants $(-)^G$ is exact, and therefore $\varphi_*(-)$ is exact.
\item From the description of $\phi_*(-)$ in (a) above, it follows that $\varphi_*(\mathcal{O}_{\left[X/G\right]}) = A^G = \mathcal{O}_{\text{Spec}(A^G)}$.
\end{enumerate}
\end{example}
Thanks to the \'etale slice theorem for algebraic stacks proved by Alper, Hall and Rydh \cite{alper-hall-rydh-etale-local}, algebraic stacks that admit good moduli spaces are \'etale locally of the form $\left[\text{Spec}(A)/G\right]$ as in the example above. Therefore the example serves as an excellent local picture to keep in mind when dealing with good moduli spaces.

The other familiar example of good moduli space comes from GIT.
\begin{example}[GIT, {\cite[\S13.5]{alper-good-moduli}}] \label{example: GIT moduli space}
Let $G$ be a linearly reductive group over a field $k$. Let $X$ be a projective-over-affine scheme over $k$ equipped with an action of $G$. Let $\mathcal{L}$ be an ample $G$-equivariant line bundle on $X$. Let $X^{\mathcal{L}\dash \rm{ss}} \subset X$ denote the open $G$-subscheme of $\mathcal{L}$-semistable points. Then the natural morphism $\phi: [X^{\mathcal{L}\dash \rm{ss}}/G] \to \text{Proj}\left(\bigoplus_{n\geq 0} H^0(X,\mathcal{L}^{\otimes n})^G\right)$ is a good moduli space for the stack $[X^{\mathcal{L}\dash \rm{ss}}/G]$.
\end{example}

Note that \Cref{example: good moduli space for affine over reductive} is a special case of \Cref{example: GIT moduli space}, where we take the equivariant line bundle to be the structure sheaf $\cO_X$, so that all points are semistable.

We should stress that not all stacks admit good moduli spaces. Having a good moduli space imposes some necessary conditions on the stack $\mathcal{M}$. For example, it forces every closed point of $\mathcal{M}$ to have linearly reductive stabilizers. This leads to an important question: what properties of an algebraic stack guarantee that it admits a good moduli space? In the next subsection we shall provide an answer to this question by explaining a set of intrinsic conditions  \cite{alper2019existence} that are necessary and sufficient for a stack over $k$ to admit a separated good moduli space.
\end{subsection}
\begin{subsection}{Intrinsic construction of moduli spaces}
One of the most useful foundational results in the theory of algebraic spaces and algebraic stacks is Artin's theorem \cite{artin-criteria} \cite[\href{https://stacks.math.columbia.edu/tag/07SZ}{Tag 07SZ}]{stacks-project}. This theorem provides a set of intrinsic local criteria (now known as the Artin criteria) which guarantee that a given functor or pseudofunctor is represented by an algebraic space or an algebraic stack. Artin's criteria are one of the main tools in the construction of an algebraic stack that parametrizes a given moduli problem.

If we want to ``rigidify" the moduli problem and look for the construction of good moduli spaces, then classically the main tool available is the theorem of Keel and Mori in the case when the moduli problem has finite inertia. One can first apply Artin's criteria to show representability by an algebraic stack, and then use the Keel-Mori theorem to construct the coarse space.

In the case when the moduli problem has infinite automorphisms, a recent breakthrough paper of Alper, Halpern-Leistner and Heinloth \cite{alper2019existence} provides a set of intrinsic local conditions for an algebraic stack to admit a separated good moduli space.

The main strategy for the construction of such a space is conceptually simple. First, the \'etale slice theorem for stacks due to Alper, Hall and Rydh \cite{alper-hall-rydh-etale-slice, alper-hall-rydh-etale-local} states that, under certain hypotheses, one can present an algebraic stack $\mathcal{M}$ \'etale locally as a quotient stack of the form $\left[\text{Spec}(A)/G\right]$, where $G$ is a linearly reductive group. We know how to construct moduli spaces for such stacks, by Example \ref{example: good moduli space for affine over reductive}. Then one would hope to be able to glue such moduli spaces in order to obtain a good moduli space for the stack $\mathcal{M}$. The insight of \cite{alper2019existence} is that there are some intrinsic properties of the stack $\mathcal{M}$ that guarantee that: (1) the hypotheses of the slice theorem are satisfied; and (2) one can glue the good moduli spaces for the local models in order to obtain a separated good moduli space for $\mathcal{M}$.

This provides an intrinsic way of constructing moduli spaces associated to moduli problems in algebraic geometry. One can first apply the Artin criteria to prove that the natural groupoid-valued pseudofunctor is represented by an algebraic stack. Then, using the criteria for the existence of good moduli spaces in \cite{alper2019existence}, one can construct the sought-after moduli space.

This strategy has been useful in cases where the techniques of GIT are not available. One such example is the construction of the moduli of Bridgeland semistable objects in the derived category of a smooth projective variety \cite[\S 7]{alper2019existence}, where the stack $\mathcal{M}$ is not known to be a quotient stack. Another example is the moduli of $K$-polystable Fano varieties \cite{alper-blum-halpern-leistner-xu,openness-k-semistable}, where the line bundle inducing semistability is not ample on the natural parameter space given by the Hilbert scheme. More recently, the same techniques have also been applied to obtain intrinsic constructions of moduli of quiver representations \cite{sveta-quivers}. We would also like to mention the construction of moduli spaces of objects in abelian categories, treated in \cite{alper2019existence}. Another important question is the projectivity of the moduli spaces obtained via intrinsic methods \cite{tajakka-thesis, projectivity-moduli-vb-curves, sveta-quivers}; we will not treat this aspect of the theory in this survey. 

To end this section, we would like to describe more carefully the criteria in \cite{alper2019existence}. They consist of filling conditions with respect to some morphisms from test stacks, reminiscent of the valuative criterion for properness in classical algebraic geometry.

We start by introducing the relevant test stacks. Let $\mathbb{G}_m$ denote the multiplicative group, whose group of $S$-points is the group of invertible functions $H^0(\cO_S^{\times})$ for any $k$-scheme $S$. We let $\mathbb{G}_m$ act on the affine space $\mathbb{A}^1_k = \text{Spec}(k[t])$ linearly such that the weight of the coordinate $t$ is $-1$. 

\begin{definition}[stack $\Theta$]
\label{defn: theta-stack}
We denote by $\Theta$ the quotient stack $\left[ \mathbb{A}^1_k / \mathbb{G}_m\right]$.
\end{definition}

Notice that the stack $\Theta$ has two $k$-points: the closed point $0$ with stabilizer $\mathbb{G}_m$ corresponding to the origin in $\mathbb{A}^1_k$, and the open point $1$ with trivial stabilizer corresponding the open orbit in $\mathbb{A}^1_k$. If $R$ is a discrete valuation ring over $k$ with uniformizer $\varpi$, then the base-change $\Theta_{R} := \Theta \times_{\text{Spec(k)}} \text{Spec}(R)$ has a unique closed point $(0,0)$ corresponding to the simultaneous vanishing locus of $t$ and $\varpi$ inside $\mathbb{A}^1_{R} = \text{Spec}(R[t])$.

\begin{definition}[$\Theta$-reductivity]
\label{defn: theta reductivity}
An algebraic stack $\mathcal{M}$ over $k$ is said to be $\Theta$-reductive if for all discrete valuation rings $R$ over $k$ and all morphisms \linebreak $f: \Theta_{R} \setminus (0,0) \to \mathcal{M}$, there is a unique dashed arrow as below making the diagram commutative
\begin{figure}[H]
\centering
\begin{tikzcd}
  \Theta_R \setminus (0,0) \ar[rd, "f"] \ar[d, symbol= \hookrightarrow] & \\ \Theta_R \ar[r, dashrightarrow] & \mathcal{M}
\end{tikzcd}
\end{figure}
\end{definition}

We continue with a choice of discrete valuation ring $R$ over $k$ with uniformizer $\varpi$.
\begin{definition}[stack $\overline{ST}_{R}$]
\label{defn: STR-stack}
 We define the quotient stack 
 \[\overline{ST}_{R}:= \left[\text{Spec}(R[t,s]/(st-\varpi))/ (\mathbb{G}_m)_{R} \right],\] where $(\mathbb{G}_m)_{R} = \mathbb{G}_m \times_{\text{Spec}(k)} \text{Spec}(R)$ acts with weight $-1$ on $t$ and with weight $1$ on $s$.
\end{definition}
The stack $\overline{ST}_{R}$ has a unique closed point$(0,0)$ corresponding to the simultaneuous vanishing locus of $s,t,\varpi$ in $\text{Spec}(R[t,s]/(st-\varpi))$.
\begin{definition}[$S$-completeness]
\label{defn: S-completeness}
An algebraic stack $\mathcal{M}$ over $k$ is defined to be $S$-complete if for all discrete valuation rings $R$ over $k$ and all morphisms $f: \overline{ST}_{R} \setminus (0,0) \to \mathcal{M}$, there is a unique dashed arrow as below making the diagram commutative
\begin{figure}[H]
\centering
\begin{tikzcd}
  \overline{ST}_{R} \setminus (0,0) \ar[rd, "f"] \ar[d, symbol= \hookrightarrow] & \\ \overline{ST}_{R} \ar[r, dashrightarrow] & \mathcal{M}
\end{tikzcd}
\end{figure}
\end{definition}

We refer the reader to \cite[\S3.5]{alper2019existence} for some intuition behind the filling condition in the definition above. We are now ready to state one of the main theorems in \cite{alper2019existence}.
\begin{theorem}[{\cite[Thm. A]{alper2019existence}}] \label{thm: existence algebraic spaces}
Let $\mathcal{M}$ be an algebraic stack with affine diagonal and of finite type over a field $k$ of characteristic $0$. Then $\mathcal{M}$ admits a separated good moduli space if and only if it is $\Theta$-reductive and $S$-complete.
\end{theorem}

\begin{remark}
The pursue of the right analogue of Theorem \ref{thm: existence algebraic spaces} in positive or mixed characteristic is a question that is being actively investigated.
\end{remark}

Our main concern now is to be able to identify moduli stacks that are $\Theta$-reductive and $S$-complete. Unfortunately, many algebraic stacks arising in moduli theory (such as the moduli stack of vector bundles on a curve) are not $\Theta$-reductive or $S$-complete. This is to be expected; we know from the work of Mumford that in order to construct a moduli space for vector bundles it is necessary to restrict to those vector bundles that are (semi)stable. This phenomenon is fairly common in moduli theory. It is often impossible to parametrize all objects in the moduli problem; instead one generally restricts their attention to a subset of objects satisfying a stability condition, which are ``rigid" enough to admit a parametrization in terms of a moduli space.

Hence we put ourselves in the following situation: we start with a ``large" stack $\mathcal{M}$ parametrizing all objects for a given moduli problem. The main question is: can we identify a reasonable open substack $\mathcal{M}^{\rm{ss}} \subset \mathcal{M}$ of ``semistable objects" such that $\mathcal{M}^{\rm{ss}}$ admits a good moduli space? In the next subsection we provide an answer to this question via the notion of $\Theta$-stability.
\end{subsection}
\begin{subsection}{$\Theta$-stability}
In the setting of GIT, Mumford introduced a notion of semistability depending on the choice of an equivariant line bundle. Suppose that $X$ is a projective variety over $k$ equipped with an action of a reductive group $G$. Given a $G$-equivariant ample line bundle $\mathcal{L}$ on $X$, Mumford defined a $G$-stable open subset $X^{\rm{ss}} \subset X$ called the semistable locus. He showed that the stack $\left[X^{\rm{ss}}/G\right]$ admits a projective good moduli space. Hence GIT provides a powerful tool to identify appropriate open substacks $\left[X^{\rm{ss}}/G\right] \subset \left[X/G\right]$ depending on an equivariant line bundle $\mathcal{L} \in \Pic^{G}(X)$, which corresponds to a line bundle on the stack $\left[X/G\right]$.

The notion of $\Theta$-stability, introduced by Halpern-Leistner \cite{halpernleistner2021structure} and Heinloth \cite{heinloth2018hilbertmumford}, generalizes the notion of GIT stability to the setting of arbitrary algebraic stacks. It is inspired by the Hilbert-Mumford criterion \cite[Chpt. 2.1]{mumford-git}, which describes the semistable locus in GIT in terms of one-parameter subgroups of the reductive group $G$. In the case of a more general algebraic stack, the one-parameter subgroups are replaced by the notion of a $\Theta$-filtration. This uses the test stack $\Theta$ introduced in Definition \ref{defn: theta-stack}.

\begin{definition}[$\Theta$-Filtration]
\label{def:theta-filtration}
Let $\mathcal{M}$ be an algebraic stack over $k$. Let $K \supset k$ be a field over $k$, and let $p: \text{Spec}(K) \to \mathcal{M}$ be a $K$-point in $\mathcal{M}(K)$. A $\Theta$-filtration of $p$ consists of the data of a morphism $f: \Theta_K \to \mathcal{M}$ along with an isomorphism $\Theta|_{1_K} \simeq p$.
\end{definition}

\begin{remark}
In \cite{halpernleistner2021structure} a $\Theta$-filtration is simply called a filtration. We stress the name $\Theta$-filtration to clearly distinguish this intrinsic notion of filtration from the many other filtrations appearing in this survey.
\end{remark}

Fix a line bundle $\mathcal{L}$ on the stack $\mathcal{M}$. For any filtration $f: \Theta_K \to \mathcal{M}$ of a point $p \in \mathcal{M}(K)$, the $0_K$-fiber of the pullback $f^*(\mathcal{L})$ acquires a natural action of the stabilizer $(\mathbb{G}_{m})_K$. We denote by $\text{wt}\left(f^*(\mathcal{L})|_{0_K}\right) \in \mathbb{Z}$ the weight of the $\mathbb{G}_m$ action on the one dimensional $0_K$-fiber.

\begin{definition}[$\Theta$-stability]
\label{def:theta-stability}
Let $\mathcal{M}$ be an algebraic stack over $k$, and let $\mathcal{L}$ be a line bundle on $\mathcal{M}$. Fix an algebraically closed field $K \supset k$ and a geometric point $p: \text{Spec}(K) \to \mathcal{M}$. We say that $p$ is $\Theta$-semistable (with respect to $\mathcal{L}$) if for all $\Theta$-filtrations $f: \Theta_K \to \mathcal{M}$ of the point $p$ we have $\text{wt}\left(f^*(\mathcal{L})|_{0_K}\right) \leq 0$.
\end{definition}

The set of $\Theta$-semistable geometric points with respect to a line bundle $\mathcal{L}$ yields a subset of the topological space of the stack $|\mathcal{M}|$ (as in \cite[\href{https://stacks.math.columbia.edu/tag/0DQM}{Tag 0DQM}]{stacks-project}). We say that the $\Theta$-semistable locus is open if this subset is open. In this case the $\Theta$-semistable points consist exactly of the geometric points of an open substack, which we denote by $\mathcal{M}^{\mathcal{L}\dash \rm{ss}} \subset \mathcal{M}$ (or just by $\mathcal{M}^{\rm{ss}}$, if the choice of the line bundle is implicit).

\begin{example}[Stability for vector bundles, \cite{heinloth2018hilbertmumford}] Let $\Bun_{n}(C)$ denote the stack of rank $n$ vector bundles on a smooth projective curve $C$ over $k$. There is a line bundle $L_{det} \in \Pic(\Bun_{n}(C))$ called the determinant line bundle, defined as follows.

By definition, there is a universal rank $n$ vector bundle $\mathcal{E}_{univ}$ on the product $\Bun_{n}(C) \times C$. The derived push-forward $R(\text{pr}_1)_*\mathcal{E}_{univ}$ under the first projection $\text{pr}_1: \Bun_{n}(C) \times C \to \Bun_{n}(C)$ is a perfect complex on $\Bun_{n}(C)$, since the morphism $\text{pr}_1$ is proper, schematic and flat. We set $L_{det}$ to be the determinant $L_{det} = \text{det}(R(\text{pr}_1)_*\mathcal{E}_{univ})$ of the perfect complex $R(\text{pr}_1)_*\mathcal{E}_{univ}$.

In \cite{heinloth2018hilbertmumford} it is proven that the semistable locus in $\Bun_{n}(C)$ with respect to the line bundle $L_{det}$ is open and consists of the slope semistable vector bundles in the sense of Mumford. We already know by GIT that this substack of semistable vector bundles $\Bun_{n}(C)^{\rm{ss}}$ admits a good moduli space.
\end{example}

Now equipped with the notion of $\Theta$-semistability, we have new tools to identify open substacks inside $\mathcal{M}$ that have a chance at admitting good moduli spaces. The question becomes: how can we check the sufficient existence criteria ($\Theta$-reductivity and $S$-completeness) for the stack of $\Theta$-semistable points $\mathcal{M}^{\rm{ss}} \subset \mathcal{M}$? We will answer this question later in this survey; this is Halpern-Leistner's Theorem \ref{thm: theta stability paper theorem} for intrinsic GIT. In order to state this theorem, we will need to gain some understanding of the complement of the $\Theta$-semistable locus $\mathcal{M} \setminus \mathcal{M}^{\rm{ss}}$, which is called the $\Theta$-unstable locus. And to acquire knowledge about the unstable locus and, as a consequence, gain insight on the $\Theta$-semistable locus, Halpern-Leistner \cite{halpernleistner2021structure} introduced the concept of $\Theta$-stratifications, to be explained in the next section.
\end{subsection}
\begin{subsection}{$\Theta$-stratifications}
\label{ssec: theta stratifications}
In this section we explain the notions of $\Theta$-stratifications and numerical invariants introduced \cite{halpernleistner2021structure}. This aims to generalize to the setting of stacks the stratifications from classical GIT \cite{shatz,kempf-filtration, hesselink-concentration, ness-stratification, kirwan-cohomology-quotients, nitsurehnsheaves}, which have also been translated to non-algebraic gauge theoretical frameworks \cite{Bruasse2005}. For simplicity we work with an algebraic stack $\mathcal{M}$ with affine diagonal and locally of finite type over a field $k$ of characteristic $0$.

We denote by $\text{Filt}(\mathcal{M}) = \text{Map}(\Theta, \, \mathcal{M})$ the stack of $\Theta$-filtrations. It is represented by an algebraic stack locally of finite type over $k$ \cite[Prop. 1.1.2]{halpernleistner2021structure}. There is a morphism $\text{ev}_1: \text{Filt}(\mathcal{M}) \to \mathcal{M}$ given by evaluating at $1 \in \Theta$.

% A graded point of the stack $\mathcal{M}$ is a morphism $(B\mathbb{G}_m)_{K} \rightarrow \mathcal{M}$, where $K \supset k$ is a field. The mapping stack of graded points $\text{Grad}(\mathcal{M}) \vcentcolon = \text{Map}(B\mathbb{G}_m, \, \mathcal{M})$ is also an algebraic stack locally of finite type over $k$ \cite[Prop. 1.1.2]{halpernleistner2021structure}.

\begin{definition}[$\Theta$-stratum, {\cite[Defn. 2.1.1]{halpernleistner2021structure}}]
Let $\mathcal{X}$ be an open substack of $\mathcal{M}$. A $\Theta$-stratum of $\mathcal{X}$ is a union of connected components of $\text{Filt}(\mathcal{X})$ such that the restriction of $\text{ev}_1$ is a closed immersion.
\end{definition}

We can think of a $\Theta$-stratum as a closed substack of $\mathcal{X}$ that is identified with some connected components of $\text{Filt}(\mathcal{X})$.

\begin{definition}[$\Theta$-stratification, {\cite[Defn. 2.1.2]{halpernleistner2021structure}}] \label{defn: theta stratification}
A $\Theta$-stratification of $\mathcal{M}$ consists of a collection of open substacks $(\mathcal{M}_{\leq c})_{c \in \Gamma}$ indexed by a totally ordered set $\Gamma$. We require the following conditions to be satisfied
\begin{enumerate}[(1)]
\item $\mathcal{M}_{\leq c} \subset \mathcal{M}_{\leq c'}$ for all $c< c'$.
\item $\mathcal{M} = \bigcup_{c \in \Gamma} \mathcal{M}_{\leq c}$.
\item For all $c$, there exists a $\Theta$-stratum $\mathfrak{S}_c \subset \text{Filt}(\mathcal{M}_{\leq c})$ of $\mathcal{M}_{\leq c}$ such that
\[ \mathcal{M}_{\leq c} \setminus \text{ev}_1(\mathfrak{S}_c) = \bigcup_{c' < c} \mathcal{M}_{\leq c'}.\]
\item For every point $p \in \mathcal{M}$, the set $\left\{ c \in \Gamma \, \mid \, p \in \mathcal{M}_{\leq c}\right\}$ has a minimal element.
\end{enumerate}
\end{definition}

A useful way to define $\Theta$-stratifications is via the use of numerical invariants in the sense of \cite[Defn. 0.0.3]{halpernleistner2021structure}. Let us give in this paragraph a simplified account of how this works for intuition; a more precise mathematical discussion will follow afterwards. A numerical invariant yields a locally constant function $\nu$ on the topological space $|\text{Filt}(\mathcal{M})|$ of the stack of $\Theta$-filtrations. This function is taken to be valued in a totally ordered group $\Gamma$ (usually $\Gamma = \mathbb{R}$ or more generally $\Gamma = \mathbb{R}[n]$ for polynomial numerical invariants). For any given geometric point $p$ of the stack $\mathcal{M}$, we can consider the subspace $|\text{Filt}(\mathcal{M}, p)| \subset |\text{Filt}(\mathcal{M})|$ consisting of the images of all $\Theta$-filtrations of $p$ (by forgetting the identification $f(1) \simeq p$). It is often the case that the function $\nu$ attains a maximum on $|\text{Filt}(\mathcal{M}, p)|$, which roughly measures the ``instability" of the point $p$. For each $c \in \Gamma$, we can then define $|\mathcal{M}_{\leq c}| \subset |\mathcal{M}|$ to be the set of geometric points $p$ such that the maximum of $\nu$ on $|\text{Filt}(\mathcal{M}, p)|$ (i.e. the ``degree of instability") is at most $c$. Under favorable circumstances, the subspace $|\mathcal{M}_{\leq c}| \subset |\mathcal{M}|$ is actually the set of geometric points of an open substack $\mathcal{M}_{\leq c} \subset \mathcal{M}$, and we obtain a $\Theta$-stratification of $\mathcal{M}$. 

For the rest of this section we shall give a more mathematically precise account of the previous paragraph. We shall restrict ourselves with numerical invariants valued in the group $\Gamma = \mathbb{R}[n]$ of polynomials in the variable $n$ with coefficients in $\mathbb{R}$. This is sufficient for many applications, such as in the case of Gieseker stability discussed in later sections.

\begin{definition}[Total order on polynomials] 
\label{defn: order_on_polynomials}
For any two polynomials $p(n),q(n) \in \mathbb{R}[n]$, we write $p(n) \leq q(n)$ if there exists some $N$ such that $p(m) \leq q(m)$ for all $m \in \mathbb{R}$ with $m>N$.
\end{definition}

The relation ``$\leq$" in Definition \ref{defn: order_on_polynomials} equips $\mathbb{R}[n]$ with the structure of a totally ordered abelian group.

\begin{definition}[Non-degenerate graded points]
\label{defn: graded_points}
For $q\geq 1$, let \linebreak $B\mathbb{G}^q_m \vcentcolon = \left[\Spec(k) / \mathbb{G}^q_m \right]$. Let $K \supset k$ be a field. We say that a morphism \linebreak $g: B(\mathbb{G}^q_m)_{K} \rightarrow \mathcal{M}$ is non-degenerate if the induced homomorphism \linebreak $\gamma: (\mathbb{G}_{m}^q)_{K} \rightarrow \text{Aut}(g(\Spec(K)))$ has finite kernel.
\end{definition}

We note that the data of a morphism $B(\mathbb{G}^q_m)_{K} \to \cM$ as in the definition above is equivalent to a pair $(p, \gamma)$ of a point $p \in \cM(K)$ and a homomorphism $\gamma: (\mathbb{G}_m^q)_{K} \to \text{Aut}(p)$.

A numerical invariant is an assignment of certain scale-invariant functions to each non-degenerate graded point $g: B(\mathbb{G}^q_m)_{K} \rightarrow \mathcal{M}$. A more formal definition follows.

\begin{definition}[Polynomial numerical invariant] \label{defn: numerical poly invariant}
A polynomial numerical invariant $\nu$ on the stack $\mathcal{M}$ is an assignment defined as follows. Let $K \supset k$ be a field and let $p \in \mathcal{M}(K)$. Let $\gamma: (\mathbb{G}_m^q)_{K} \rightarrow \text{Aut}(p)$ be a homomorphism with finite kernel. The polynomial numerical invariant $\nu$ assigns to this data a scale-invariant function $\nu_{\gamma}: \mathbb{R}^q\setminus \{0\} \rightarrow \mathbb{R}[n]$ such that:
\begin{enumerate}[(1)]
    \item $\nu_{\gamma}$ is unchanged under field extensions $K \subset K'$.
    \item $\nu$ is locally constant in algebraic families. In other words, choose any scheme $T$, a morphism $\xi: T \rightarrow \mathcal{M}$, and a homomorphism $\gamma: (\mathbb{G}_{m}^{q})_{T} \rightarrow \text{Aut}(\xi)$ of $T$-group schemes with finite kernel. Then as we vary $t \in T$, we require that the function $\nu_{\gamma_{t}}$ is locally constant in $T$.
    \item Given a homomorphism $\phi: (\mathbb{G}^w_{m})_K \rightarrow (\mathbb{G}^q_{m})_{K}$ with finite kernel, the function $\nu_{\gamma \circ \phi}$ is the restriction of $\nu_{\gamma}$ along the inclusion $\mathbb{R}^w \hookrightarrow \mathbb{R}^q$ induced by $\phi$.
\end{enumerate}
\end{definition}

A useful way to construct polynomial numerical invariants on a stack $\mathcal{M}$ is to use line bundles on $\mathcal{M}$, which are the analogues of equivariant polarizations in the classical setting of GIT. Fix a sequence of rational line bundles $(\mathcal{L}_n)_{n \in \mathbb{Z}}$, where each $\mathcal{L}_n$ is in the rational Picard group $\Pic(\mathcal{M})\otimes_{\mathbb{Z}} \mathbb{Q}$ of the stack. For each morphism $g: (B\mathbb{G}^q_m)_{K} \rightarrow \mathcal{M}$, the pullback line bundle $g^*(\mathcal{L}_{n})$ amounts to a rational character in $X^*(\mathbb{G}^q_m)\otimes_{\mathbb{Z}} \mathbb{Q}$. Under the natural identification $X^*(\mathbb{G}^q_m)\otimes_{\mathbb{Z}} \mathbb{Q} \cong \mathbb{Q}^q$, we can interpret this as a $q$-tuple of rational numbers $(w_n^{(i)})_{i=1}^q$, which we call the weight of $g^*(\mathcal{L}_n)$. 

In practice we can often choose the line bundles $(\mathcal{L}_n)_{n \in \mathbb{Z}}$ in such a way that, for each fixed $i$, there is a polynomial $w^{(i)}(T) \in \mathbb{R}[T]$ such that the weight $w_n^{(i)}$ is $w^{(i)}(n)$. If this is the case, then for every $g$ we can define an $\mathbb{R}$-linear function $L_{g} : \mathbb{R}^q \to \mathbb{R}[n]$ given by
\[ L_{g}\left((r_i)_{i=1}^q\right) = \sum_{i=1}^q r_i \cdot w^{(i)}(n).\]
In order to obtain a scale invariant $\nu$, we use a rational quadratic norm on graded points as in \cite[Defn. 4.1.12]{halpernleistner2021structure}, which is the analogue of a Weyl invariant quadratic form in the classical setting of GIT. 

\begin{definition}[Rational quadratic norm on graded points]
\label{defn: rational_quadratic_norm}
A rational quadratic norm $b$ on graded points of the stack $\mathcal{M}$ is an assignment defined as follows. Let $K \supset k$ be a field and let $p \in \mathcal{M}(K)$. Let $\gamma: (\mathbb{G}_m^q)_{K} \rightarrow \text{Aut}(p)$ be a homomorphism with finite kernel. Then $b$ assigns to this data a positive definite quadratic norm $b_{\gamma}: \mathbb{R}^q \to \mathbb{R}$ with rational coefficients such that:
\begin{enumerate}[(1)]
    \item $b_{\gamma}$ is unchanged under field extensions $K \subset K'$.
    \item $b$ is locally constant in algebraic families. In other words, choose any scheme $T$, a morphism $\xi: T \rightarrow \mathcal{M}$ and a homomorphism $\gamma: (\mathbb{G}_{m}^{q})_{T} \rightarrow \text{Aut}(\xi)$ of $T$-group schemes with finite kernel. Then as we vary $t \in T$, we require that the function $b_{\gamma_{t}}$ is locally constant in $T$.
    \item Given a homomorphism $\phi: (\mathbb{G}^w_{m})_K \rightarrow (\mathbb{G}^q_{m})_{K}$ with finite kernel, the norm $b_{\gamma \circ \phi}$ is the restriction of $b_{\gamma}$ along the inclusion $\mathbb{R}^w \hookrightarrow \mathbb{R}^q$ induced by $\phi$.
\end{enumerate}
\end{definition}

\begin{constr} \label{constr: numerical invariant}
Given a sequence of rational line bundles $(\mathcal{L}_n)_{n \in \mathbb{Z}}$ as above and a rational quadratic norm on graded points $b$, we define a polynomial numerical invariant $\nu$ as follows. For all non-degenerate $g: (B\mathbb{G}^q_m)_{K} \rightarrow \mathcal{M}$ with corresponding morphism $\gamma: (\mathbb{G}^q_m)_{K} \to \text{Aut}(g(\Spec(K)))$, we set
\[ \nu_{\gamma}(\vec{r}) = \frac{L_{g}(\vec{r})}{\sqrt{b_{\gamma}(\vec{r})}}\]
where $\overrightarrow{r}=(r_1,\ldots, r_q)$.
\end{constr}

We next explain how a polynomial numerical invariant $\nu$ can be used top define a $\Theta$-stratification on $\mathcal{M}$ (see \cite[\S 4.1]{halpernleistner2021structure} for more details).  We say that a $\Theta$-filtration $f: \Theta_{K} \rightarrow \mathcal{M}$ is non-degenerate if the restriction $f|_0: [0/\, (\mathbb{G}_m)_{K}] \to \mathcal{M}$ is non-degenerate. We regard $\nu$ as a function on the set of non-degenerate $\Theta$-filtrations by defining $\nu(f) \vcentcolon = \nu_{f|_{0}}(1) \in \mathbb{R}[n]$.
\begin{definition}[Semistability for numerical invariants]
\label{defn: semistability-nu}
Given a polynomial numerical invariant $\nu$ on $\mathcal{M}$ and a point $p \in |\mathcal{M}|$, we say that $p$ is $\nu$-semistable if all non-degenerate $\Theta$-filtrations $f$ with $f(1)\simeq p$ satisfy $\nu(f) \leq 0$. Otherwise we say that $p$ is $\nu$-unstable. 
\end{definition}

For any $\nu$-unstable point $p$, we set $M^{\nu}(p)$ to be the supremum of $\nu(f)$ over all $\Theta$-filtrations $f$ with $f(1)\simeq p$ (if such supremum exists). If $p$ is $\nu$-semistable, then by convention we set $M^{\nu}(p) = 0$. For any $c \in \mathbb{R}[n]_{\geq 0}$, let $\mathcal{M}_{\leq c}$ be the set of all points $p$ satisfying $M^{\nu}(p) \leq c$, and we let $\cM^{\nu\dash \rm{ss}} := \cM_{\leq 0}$ denote the set of $\nu$-semistable points. If this is an open stratification of the stack coming from a $\Theta$-stratification, then we say that $\nu$ defines a $\Theta$-stratification.

\begin{remark}
\label{remark: non-deg filtrations}
Generally one only considers non-degenerate $\Theta$-filtrations, because these are the ones relevant for stability. Hence the adjective ``non-degenerate" is usually omitted.
\end{remark}

\begin{definition}[$\nu$-Harder-Narasimhan filtration] \label{defn: hn-filtration}
For any $\nu$-unstable point $p$, a $\Theta$-filtration $f$ of $p$ is called a $\nu$-Harder-Narasimhan filtration, abbreviated $\nu$-HN filtration, if $\nu(f) = M^{\nu}(p)$. 
\end{definition}

If $\nu$ defines a $\Theta$-stratification, then by \cite[Lemma 2.1.4]{halpernleistner2021structure} every $\nu$-unstable point admits a $\nu$-Harder-Narasimhan filtration that is unique up to pre-composing with a ramified covering $\Theta \to \Theta$. 

We should remark that not all numerical invariants $\nu$ coming from sequences of line bundles and rational quadratic norms define $\Theta$-stratifications. In the next section we discuss some criteria that guarantee the existence of such a $\Theta$-stratification and, in addition, imply that the open semistable locus $\mathcal{M}^{\nu\dash \rm{ss}}$ admits a good moduli space.
\end{subsection}

\begin{subsection}{Monotonicity and intrisic moduli theory}
In this section we explain Halpern-Leistner's intrinsic GIT theory \cite[\S 5]{halpernleistner2021structure}. We keep our assumption that $\mathcal{M}$ is an algebraic stack with affine diagonal and locally of finite type over a field $k$ of characteristic $0$.

 A natural question to ask is: when does a polynomial numerical invariant $\nu$ define a $\Theta$-stratification as described in the previous section? Another thing we would like to know is whether the locus $\mathcal{M}^{\nu \dash \rm{ss}}$ of $\nu$-semistable objects admits a good moduli space. A general strategy to address such questions proceeds in three steps:
 \begin{enumerate}[(Step 1)]
    \item Show that the polynomial numerical invariant $\nu$ satisfies some monotonicity conditions (as in Definition \ref{defn: strictly theta monotone and STR monotone}).
     \item Prove that $\nu$ satisfies the HN boundedness condition (as in Definition \ref{defn: HN boundedness}).
    \item Check that the semistable locus $\mathcal{M}^{\nu \dash \rm{ss}}$ is a disjoint union of quasi-compact open substacks.
 \end{enumerate}

The following theorem describes this strategy more precisely.
\begin{theorem}[Intrinsic GIT, {\cite[Theorem B]{halpernleistner2021structure}}] \label{thm: theta stability paper theorem}
Let $\mathcal{M}$ be an algebraic stack with affine diagonal and locally of finite type over a field $k$ of characteristic $0$. Let $\nu$ be a polynomial numerical invariant on $\mathcal{M}$ defined by a sequence of rational line bundles and a norm on graded points, as explained in the Construction \ref{constr: numerical invariant}. Then,
\begin{enumerate}
    \item If $\nu$ is strictly $\Theta$-monotone (Definition \ref{defn: strictly theta monotone and STR monotone}) and satisfies the HN boundedness condition (Definition \ref{defn: HN boundedness}), then it defines a $\Theta$-stratification of $\mathcal{M}$.
    \item Assume that $\nu$ satisfies the conditions in (1). In addition, suppose that $\nu$ is strictly $S$-monotone (Definition \ref{defn: strictly theta monotone and STR monotone}) and that the semistable locus $\mathcal{M}^{\nu \dash \rm{ss}}$ can be written as a disjoint union of quasi-compact open substacks. Then $\mathcal{M}^{\nu \dash \rm{ss}}$ has a separated good moduli space.
    \item If all the conditions in (2) are satisfied and if $\mathcal{M}$ satisfies the existence part of the valuative criterion for properness for complete discrete valuation $k$-algebras \cite[\href{https://stacks.math.columbia.edu/tag/0CLK}{Tag 0CLK}]{stacks-project}, then each quasi-compact open and closed substack of $\mathcal{M}^{\nu \dash \rm{ss}}$ admits a good moduli space that is proper over $k$.
\end{enumerate}
\end{theorem}

 Recall that, in Mumford's GIT, the quotient of a projective variety with an ample linearization is automatically projective. In many applications of GIT, this is the way in which properness of the moduli space is proved.  But in Theorem \ref{thm: theta stability paper theorem} we see that, in order to prove properness, we have to show a valuative criterion. The prototype of such criteria is  Langton's semistable reduction theorem for torsion free sheaves \cite{langton}. There has been a lot of work extending this to principal bundles in different contexts (cf. \cite{balaji.seshadri.1,balaji.parameswaran.2,balaji.du,heinloth.1,heinloth.2}).

We end this subsection by explaining each of the hypotheses that one needs to check in Theorem \ref{thm: theta stability paper theorem}. We need some setup first. 

\begin{notation}[``Rigidified" $\Theta$ and $\overline{ST}_{R}$] \label{notn: rigidified theta and str}
Let $R$ be a complete discrete valuation ring over $k$. Choose a uniformizer $\varpi$ of $R$. We define $Y_{\Theta_{R}} \vcentcolon = \mathbb{A}^1_{R}$ equipped with the $\mathbb{G}_m$-action that gives $t$ weight $-1$. We have $\Theta_{R} = \left[ \, \mathbb{A}^1_{R} / \, \mathbb{G}_m \, \right]$. Note that $\mathbb{A}^1_{R}$ contains a unique $\mathbb{G}_m$-invariant closed point cut out by the ideal $(t, \varpi)$. We will denote this point by $(0,0)$.

 We define $Y_{\overline{ST}_{R}} \vcentcolon = \Spec\left( \, R[t,s]/(st-\varpi) \,\right)$. We equip $Y_{\overline{ST}_{R}}$ with the $\mathbb{G}_m$-action that gives $s$ weight $1$ and $t$ weight $-1$. The isomorphism class of the $\mathbb{G}_m$-scheme $Y_{\overline{ST}_{R}}$ is independent of the choice of uniformizer $\varpi$. Define $\overline{ST}_{R} \vcentcolon = [\,Y_{\overline{ST}_{R}} / \, \mathbb{G}_m \, ]$. We will denote by $(0,0)$ the unique $\mathbb{G}_m$-fixed closed point in $Y_{\overline{ST}_{R}}$. We have that $(0,0)$ is cut out by the ideal $(s,t)$.
\end{notation}

\begin{notation}[Weighted lines]
Let $\kappa \supset k$ be a field and let $a \geq 1$ be an integer. We denote by $\mathbb{P}^1_{\kappa}[a]$ the $\mathbb{G}_m$-scheme $\mathbb{P}^{1}_{\kappa}$ equipped with the $\mathbb{G}_m$-action determined by the equation $t \cdot [x:y] = [t^{-a}x : y]$. We set $0 = [0: 1]$ and $\infty = [1:0]$.
\end{notation}

In the following definition we describe a simplified version of the monotonicity conditions in \cite[\S 5]{halpernleistner2021structure}.

\begin{definition}[Monotonicity] \label{defn: strictly theta monotone and STR monotone} A polynomial numerical invariant $\nu$ on $\mathcal{M}$ is strictly $\Theta$-monotone (resp. strictly $S$-monotone) if the following condition holds. 

Let $R$ be any complete discrete valuation ring and set $\mathfrak{X}$ to be $\Theta_{R}$ (resp. $\overline{ST}_{R}$), and $Y_{\mathfrak{X}}$ as in Notation \ref{notn: rigidified theta and str}. Choose a map $\varphi: \mathfrak{X} \setminus (0,0) \rightarrow \mathcal{M}$. Then, after maybe replacing $R$ with a finite DVR extension, there exists a reduced and irreducible scheme $\Sigma$ equipped with a $\mathbb{G}_m$-action along with $\mathbb{G}_m$-equivariant maps $f: \Sigma \rightarrow Y_{\mathfrak{X}}$ and $\widetilde{\varphi}: \left[\Sigma/ \, \mathbb{G}_m\right] \rightarrow  \mathcal{M}$ such that
\begin{enumerate}[({M}1)]
    \item The map $f$ is proper, $\mathbb{G}_m$-equivariant, and its restriction induces an isomorphism $f : \, \Sigma_{Y_{\mathfrak{X}} \setminus (0,0)} \xrightarrow{\sim} Y_{\mathfrak{X}} \setminus (0,0)$.
    \item The following diagram commutes
\begin{figure}[H]
\centering
\begin{tikzcd}
  \left[\left(\Sigma_{Y_{\mathfrak{X}} \setminus (0,0)}\right)/ \, \mathbb{G}_m \right] \ar[rd, "\widetilde{\varphi}"] \ar[d, "f"'] & \\   \mathfrak{X} \setminus (0,0) \ar[r, "\varphi"'] &  \mathcal{M}
\end{tikzcd}
\end{figure}
    \item Let $\kappa$ denote a finite extension of the residue field of $R$. For any $a \geq 1$ and any finite $\mathbb{G}_m$-equivariant morphism $\mathbb{P}^1_{\kappa}[a] \to \Sigma_{(0,0)}$, we have $\nu\left( \;\widetilde{\varphi}|_{\left[\infty / \mathbb{G}_m\right]} \;\right) >  \nu\left(\; \widetilde{\varphi}|_{\left[0 / \mathbb{G}_m\right]} \;\right)$.
\end{enumerate}
\end{definition}

\begin{definition}[HN Boundedness] \label{defn: HN boundedness}
We say that a polynomial numerical invariant $\nu$ satisfies the HN boundedness condition if the following is always satisfied:

(HNB) Let $T$ be an affine Noetherian scheme, equipped with morphism $g: T \rightarrow \mathcal{M}$. Then there exists a quasi-compact open substack $\mathcal{U}_{T} \subset \mathcal{M}$ such that the following holds.
For all geometric points $t \in T$ with residue field $\kappa(t)$ and all $\Theta$-filtrations \linebreak $f: \Theta_{\kappa(t)} \rightarrow \mathcal{M}$ of the point $g(t)$ with $\nu(f)>0$, there exists another filtration $f'$ of $g(t)$ satisfying $\nu(f') \geq \nu(f)$ and $f'|_{0} \in \mathcal{U}_{T}$.
\end{definition}
This says that for the purposes of maximizing $\nu(f)$ among all $\Theta$-filtrations of points in a bounded family (i.e. finding $\nu$-HN filtrations as in Definition \ref{defn: hn-filtration}), it suffices to consider only those $\Theta$-filtrations $f$ such that the associated graded $f|_{0}$ lies in some other (possibly larger) bounded family.
\end{subsection}
\end{section}

\begin{section}{Infinite dimensional GIT for torsion-free sheaves}
\label{sec: torsion free}
In this section we illustrate the previous notions with one of the main examples in moduli theory: the moduli problem of torsion-free sheaves. This moduli space was constructed using GIT by Mumford \cite{MumfordICM} (stable vector bundles) and Seshadri (compactification with semistable vector bundles) \cite{seshadri}  over curves,  Gieseker \cite{gieseker_torsion_free} for surfaces, later Maruyama \cite{maruyama_moduli} and Simpson \cite{simpson-repnfundamental-I} in higher dimension. 

\begin{subsection}{Torsion-free sheaves} 
\label{ssec: tf sheaves}

Let $X$ be a smooth projective variety of dimension $d$ over the field $k$, and fix an ample line bundle $\mathcal{O}(1)$ on $X$. Let $T$ be a $k$-scheme and let $\mathcal{F}$ be a $\mathcal{O}_{X_T}$-module. We say that $\mathcal{F}$ is a relative torsion-free sheaf of rank $r$ on $X_{T}$ if it is $T$-flat, finitely presented, and for all points $t \in T$ the fiber $\mathcal{F}|_{X_t}$ is torsion-free of rank $r$ (i.e. the dimension of the generic fiber is $r$).

% \andres{I propose to move this to section 4.2. Also I would mention that the determinant is a line bundle.} We describe the notion of a determinant for a torsion-free sheaf $\mathcal{F}$ of rank $r$, recalling \cite[\S 5.6]{kobayashi-vectorbundles} and \cite[Def 2.7]{rho-sheaves}. If we denote by $\bigwedge^r \mathcal{F}$ the 
% quotient of $\mathcal{F}^r$ by the submodule generated by sections 
% \[\sigma(x_1\otimes x_2\otimes \cdots \otimes x_r) - {\rm sign}(\sigma) x_1\otimes x_2\otimes\cdots \otimes x_r,\] 
% where $\sigma$
% is an element of the symmetric group on $r$ letters, we define the determinant of $\mathcal{F}$ as ${\rm det}:=\left( \bigwedge^r \mathcal{F}\right)^{\vee \vee}$, where $\vee$ denotes the dual. 

Let us present the stack of coherent sheaves that is relevant for this moduli problem. 

\begin{definition}[{\cite[Defn. 2.9]{rho-sheaves}}] \label{defn: stack of torsion free sheaves with trivialization}
Denote by $\Coh_r^{tf}(X)$ the \linebreak pseudofunctor from $\left(\Aff_{k}\right)^{op}$ into groupoids such that, for any affine scheme $T \in \Aff_{k}$, it gives
\begin{gather*} \Coh_r^{tf}(X)\, (T) \; = \; \left[ \begin{matrix} \;  \text{groupoid of $T$-flat relative torsion-free sheaves $\mathcal{F}$ of rank $r$ on $X_{T}$}  \;  \end{matrix} \right]\end{gather*}
\end{definition}
\end{subsection}

The stack $\Coh_r^{tf}(X)$ is an open substack of $\Coh_r(X)$, the algebraic stack of coherent sheaves in $X$, and contains the open substack $\Bun_r(X)$ of rank $r$ vector bundles on $X$.

Denote by $P_{\mathcal{F}}$ the Hilbert polynomial of a torsion-free sheaf $\mathcal{F}$ on $X$ with respect to the line bundle $\mathcal{O}(1)$. This is a polynomial of degree $d$ in the variable $n$
\[ P_{\mathcal{F}}(n) = \sum_{i = 0}^{d} \frac{a_i}{i!} n^i\]
with rational coefficients $a_i$ (see \cite[Lem. 1.2.1]{huybrechts.lehn}). We define the reduced Hilbert polynomial as $\overline{p}_{\mathcal{F}} \vcentcolon = \frac{1}{a_{d}} P_{\mathcal{F}}$. 
For each $0 \leq i \leq d-1$, define the $i^{th}$ slope of $\mathcal{F}$ as $\widehat{\mu}_i(\mathcal{F}) \vcentcolon = \frac{a_{i}}{a_d}$, where $\widehat{\mu}_{d-1}(\mathcal{F}) = \frac{a_{d-1}}{a_d}$ is the usual slope of $\mathcal{F}$.

We define a family of line bundles on $\Coh_r^{tf}(X)$ as follows. 

\begin{definition}[{\cite[Defn. 4.1]{rho-sheaves}}] \label{defn: line bundle torsion-free}
Given a $k$-scheme $T$, let \linebreak $\pi_T:X_T\rightarrow T$ the structure morphism. Let $f:T\rightarrow \Coh_r^{tf}(X)$ be a morphism represented by a rank $r$ torsion-free sheaf $\mathcal{F}$ of dimension $d$ on $X_T$. For each integer $n$, define a line bundle $f^*M_n$ on $T$ by
\[f^{\ast}M_n  := \det R_{\pi_T\ast}(\mathcal{F}(n)).\]
The line bundles $f^*M_n$ for all $f:T\rightarrow \Coh_r^{tf}(X)$ are compatible with pullbacks and yield a line bundle $M_n$ on $\Coh_r^{tf}(X)$. We set $b_d := \bigotimes_{i = 0}^d M_{i}^{(-1)^{d-i} \binom{d}{i}}$, and we define 
\[L_n = M_n \otimes b_d^{-\otimes \, \overline{p}_{\mathcal{F}}(n)}.\]
\end{definition}

The following step in the construction is to define $\Theta$-filtrations of the stack of torsion-free coherent sheaves, as in Definition \ref{def:theta-filtration}. A $\Theta$-filtration \linebreak $f: \Theta_K \to \Coh_{r}^{tf}(X)$ of a $K$-point of the stack $\mathcal{F}\in \Coh_r^{tf}$ consists of 
\begin{enumerate}[(a)]
    \item A $\mathbb{G}_m$-equivariant family $\widetilde{\mathcal{F}}$ of torsion-free sheaves on $X \times \mathbb{A}^1_{K}$.
    \item An identification $\widetilde{\mathcal{F}}|_{X \times 1} \xrightarrow{\sim} \mathcal{F}$.
\end{enumerate}
By using the Rees construction \cite[Prop. 1.0.1]{halpernleistner2021structure} we can rewrite the definition of $\Theta$-filtration in a much more explicit way. 

\begin{proposition}[{\cite[Prop. 3.8]{rho-sheaves}}] \label{prop: filtrations of torsion free sheaves}
A $\Theta$-filtration $f : \Theta_{K} \rightarrow \Coh_{r}^{tf}(X)$ of $\mathcal{F}$ corresponds to a $\mathbb{Z}$-filtration $(\mathcal{F}_m)_{m \in \mathbb{Z}}$ by $\mathcal{O}_{X_K}$-subsheaves of $\mathcal{F}$ satisfying
\begin{enumerate}[(1)]
    \item $\mathcal{F}_{m+1} \subset \mathcal{F}_{m}$.
    \item $\mathcal{F}_m = 0$ for $m \gg 0$ and $\mathcal{F}_m = \mathcal{F}$ for $m \ll 0$.
    \item $\mathcal{F}_{m}/\mathcal{F}_{m+1}$ is torsion-free.
\end{enumerate}
The graded object at the origin $f|_{[0/\mathbb{G}_m]} : [0/\mathbb{G}_m] \rightarrow \Coh_{r}^{tf}(X)$ is the associated graded sheaf $\bigoplus_{m \in \mathbb{Z}} \mathcal{F}_{m} / \mathcal{F}_{m+1}$.
\end{proposition}

For each $\Theta$-filtration the restrictions $f^* L_n$ of the line bundles are given by \linebreak $\mathbb{G}_m$-equivariant line bundles on $\mathbb{A}^1_{\mathbb{K}}$ and are determined by the fiber $L_n|_0$ considered as a $1$-dimensional representation of $\mathbb{G}_m$. 
The integer classifying the character of $\mathbb{G}_m$ associated to this representation is called the weight $\text{wt}(L_n|_0)$.

On the other hand, $M_n|_0$ is isomorphic to $\det\left(R_{\pi_K\ast}(\widetilde{\mathcal{F}}|_0)(n)\right)$, where the graded sheaf is $\widetilde{\mathcal{F}}|_0=\bigoplus_{m \in \mathbb{Z}} \mathcal{F}_{m} / \mathcal{F}_{m+1}$. Given that $\mathbb{G}_m$ acts on $\mathcal{F}_m / \mathcal{F}_{m+1}$ with weight $m$, we get
\[\text{wt}(M_n|_0)=\sum_{m\in \mathbb{Z}}\left(P_{\mathcal{F}_m}(n)-P_{\mathcal{F}_{m+1}}(n)\right)\cdot m\]
and 
\[\text{wt}(b_d|_0)=\sum_{m\in \mathbb{Z}}\left(\rk_{\mathcal{F}_m}-\rk_{\mathcal{F}_{m+1}}\right)\cdot m.\]
Therefore,
\[\text{wt}(L_n|_0)=\text{wt}(M_n|_0)-\overline{p}_{\mathcal{F}}(n)\cdot \text{wt}(b_d|_0)=\]
\[\sum_{m\in \mathbb{Z}}m\cdot \left(\overline{p}_{\mathcal{F}_m / \mathcal{F}_{m+1}}(n)-\overline{p}_{\mathcal{F}}(n)\right)\cdot \rk_{\mathcal{F}_m / \mathcal{F}_{m+1}}.\]

In order to set the polynomial numerical invariant as in Definition \ref{defn: numerical poly invariant}, we define a rational quadratic norm on the graded points of $\Coh^{tf}_r(X)$ (c.f. \cite[Section 4.1]{torsion-freepaper}). For simplicity we will focus only on $B\mathbb{G}_m^q$-points with $q=1$.

If $g: (B\mathbb{G}_m)_K\rightarrow \Coh_r^{tf}(X)$ is a $\mathbb{Z}$-graded torsion-free sheaf $\overline{F}=\bigoplus_{m\in \mathbb{Z}}\overline{F}_{m}$ of rank $r$ on $X_K$, define $b_g$ to be the positive definite rational quadratic form on $\mathbb{R}$ given by 
\[b_g(v):=\sum_{m\in \mathbb{Z}^q} \rk_{\overline{\mathcal{F}}_{m}}\cdot (m\cdot v)^2,\]
which is uniquely determined by its value in $v=1$. 

\begin{definition}[{\cite[Defn. 4.3]{torsion-freepaper}}]\label{defn: numerical invariant}
Let $f: \Theta_{K} \rightarrow \Coh^{tf}_r(X)$ be a  $\Theta$-filtration corresponding to the $\mathbb{Z}$-filtration $(\mathcal{F}_m)_{m\in\mathbb{Z}}$. We define the polynomial numerical invariant $\nu(f)$ to be
 \begin{equation} \label{eqn_1}
 \nu(f) := \frac{\text{wt}(L_n|_0)}{\sqrt{b_{f|_0}(1)}}=\frac{\sum_{m \in \mathbb{Z}} m \cdot  (\overline{p}_{\mathcal{F}_m/ \mathcal{F}_{m+1}} - \overline{p}_{\mathcal{F}}) \cdot \rk(\mathcal{F}_m/\mathcal{F}_{m+1})}{\sqrt{\sum_{m \in \mathbb{Z}} m^2 \cdot \rk(\mathcal{F}_m / \mathcal{F}_{m+1})}}\, .
 \end{equation}
\end{definition}

\begin{subsection}{Affine grassmannians, monotonicity and good moduli space}

Affine grassmannians are the tool developed in \cite{torsion-freepaper} to show the monotonicity 
 properties (Definition \ref{defn: strictly theta monotone and STR monotone}) required to prove Theorem \ref{thm: theta stability paper theorem}. The idea is to build an algebraic counterpart to infinite-dimensional Atiyah-Bott-Donaldson symplectic reduction.

We can define a fibered category $\Coh^{tf}_{r, rat}(X)$ whose objects are pairs $(D,\mathcal{F})$ of a divisor and torsion-free sheaf, and a morphism is a pair
$$
(i:D_1\inj D_2, \psi:\mathcal{F}_2\to \mathcal{F}_1)
$$ 
where $i$ is an inclusion and $\psi$ is a homomorphism which is an isomorphism when restricted to the complement of $D_2$ (notice that this is not fibered in groupoids). This is called the stack of rational maps, see \cite[section 3.2]{torsion-freepaper} for more details.  We have a natural morphism $\Coh^{tf}_r(X)\rightarrow \Coh^{tf}_{r, rat}(X)$, sending a torsion free sheaf $\mathcal{F}$ to 
the pair $(\emptyset,\mathcal{F})$. The fibers of this morphism are
ind-projective ind-schemes. We call them affine Grassmannians, in analogy with the notion appearing 
in the moduli of vector bundles over curves. Changing the trivialization over the open subset, we can heuristically think of this morphism as being a local presentation of the stack $\Coh^{tf}_r(X)$ as the quotient of
infinite-dimensional spaces, the ind-schemes, by the action of
an infinite-dimensional group of rational gauge transformations.

This, together with the fact that the family of line bundles defined in the previous section is asymptotically ample on each projective stratum, will be the key to show monotonicity in \Cref{thm: theta stability paper theorem}. 
In the following, $\mathcal{F}$ will denote a family of torsion-free sheaves of rank $r$ on $X_S$, where $S$ is a quasi-compact $k$-scheme. Fix an effective Cartier divisor $D \hookrightarrow X_S$ flat over $S$ and set $Q = X_S \setminus D$.

\begin{definition}[cf. {\cite[Defn. 3.12]{torsion-freepaper}}] \label{defn: step 1 grassmannian}
We define $\Gr_{X_S, D, \mathcal{F}}$ to be the functor from $\Aff_S^{\, \, op}$ to sets such that for every affine scheme $T\in \Aff_S$, the set $\Gr_{X_S, D, \mathcal{F}} (T)$ consists of equivalence classes of pairs $(\mathcal{E}, \theta)$, where
\begin{enumerate}[(1)]
    \item $\mathcal{E}$ is a family of rank $r$ torsion-free sheaves on $X_T$.
    \item $\theta$ is morphism $\theta:  \mathcal{F}|_{X_{T}} \rightarrow  \mathcal{E}$ such that the restriction $\theta|_{Q_{T}}$ is an isomorphism.
\end{enumerate}
The pairs $(\mathcal{E}_1, \theta_1)$ and $(\mathcal{E}_2, \theta_2)$ are equivalent if there is an isomorphism $\mathcal{E}_1 \xrightarrow{\sim} \mathcal{E}_2$ that identifies $\theta_1$ and $\theta_2$.
\end{definition}

We remark that the fiber over $(D,\mathcal{F}$) of the natural morphism $\Coh^{tf}_r(X)\rightarrow \Coh^{tf}_{r, rat}(X)$ is
$\Gr_{X_S, D, \mathcal{F}}$.

Now let $j$ be the open immersion $j: Q \hookrightarrow X_{S}$ and choose $T \in \Aff_S$ and an element $(\mathcal{E}, \theta) \in \Gr_{X_S, D, \mathcal{F}}(T)$. Using that the unit $\mathcal{E} \rightarrow j_{T \, *} \, j_T^* \, \mathcal{E}$ is a monomorphism we can use the isomorphism $j_{T *}(j_{T})^*\theta$ to understand $\mathcal{E}$ as a subsheaf of $j_{T \, *} \, j_T^* \, \mathcal{F}|_{X_T}$ and have $\mathcal{F}|_{X_{T}} \subset \mathcal{E} \subset j_{T \, *} \, j_T^* \, \mathcal{F}|_{X_T}$. We use the increasing sequence of sheaves 
\[\mathcal{F}|_{X_T} \subset \mathcal{F}|_{X_T}(D_{T}) \subset \mathcal{F}_{X_T}(2D_{T}) \subset \cdots \subset \mathcal{F}|_{X_T} (nD_{T}) \subset \cdots\] 
to define the truncated subfunctor $\Gr^{\leq N}_{X, D, \mathcal{F}}\subset \Gr_{X, D, \mathcal{F}}$, where
\[ \Gr^{\leq N}_{X, D, \mathcal{F}} (T) \; := \; \left\{ \begin{matrix} \; \; \; \text{pairs $(\mathcal{E}, \theta)$ in  $\Gr_{X_S, D, \mathcal{F}}(T)$ such that} \; \; \; \\ \text{$\mathcal{E} \subset \mathcal{F}|_{X_T}(ND_T)$} \end{matrix} \right\}\]
For any rational polynomial $P \in \mathbb{Q}[n]$, define the open and closed subfunctor $\Gr_{X_S, D, \mathcal{F}}^{\leq N, P}\subset \Gr_{X_S, D, \mathcal{F}}^{\leq N}$ consisting of points $(\mathcal{E}, \theta)$ where $\mathcal{E}$ has Hilbert polynomial $P$ on every fiber.

Recall that a strict ind-scheme is a colimit of closed immersions. 

\begin{proposition}[{\cite[Prop. 3.13]{torsion-freepaper}}]\label{prop: ind-representability of step 1 grassmannian}
We have $\Gr_{X_{S}, D, \mathcal{F}} = \underset{N>0}{\colim} \; \Gr^{\leq{N}}_{X_{S}, D, \mathcal{F}}$ (as presheaves on $\Aff_S$) and, for each $N \geq 0$,
\[ \Gr_{X_S, D, \mathcal{F}}^{\leq N} = \bigsqcup_{P \in \mathbb{Q}[n]} \Gr_{X_S, D, \mathcal{F}}^{\leq N, P} \; ,\]
where each $\Gr_{X_S, D, \mathcal{F}}^{\leq N, P} $ is represented by a projective scheme of finite presentation over $S$. This induces a presentation of $\Gr_{X_{S}, D, \mathcal{F}}$ as an ind-projective strict ind-scheme over $S$.
\end{proposition}

There is a natural forgetful morphism $\Gr_{X_S, D, \mathcal{F}}\rightarrow \Coh_r^{tf}(X)$ sending $(\mathcal{E},\theta)$ to $\mathcal{E}$. We denote again by $L_n$, $b_d$ and $M_n$ the line bundles on $\Gr_{X_S, D, \mathcal{F}}$ pulled back from those on $\Coh_r^{tf}(X)$ defined before, by abuse of notation. 

\begin{proposition}[{\cite[Prop. 3.19]{torsion-freepaper}}]\label{prop: ampleness of line bundle}
There exists some integer \linebreak $m\gg0$, depending on $N \in \mathbb{Z}_{\geq 0}$ and $P \in \mathbb{Q}[n]$, such that, for all $n \geq m$, the restriction $L_n^{\vee}|_{\Gr_{X_S, D, \mathcal{F}, \sigma}^{\leq N, P}}$ is $S$-ample.
\end{proposition}

Using the affine grassmannians above and the rational filling condition proven in \cite[Lem. 4.7]{torsion-freepaper}, one can show that the stack $\Coh^{tf}_r(X)$ satisfies the two monotonicity conditions in Definition \ref{defn: strictly theta monotone and STR monotone}. This type of argument was called ``Infinite dimensional GIT" in \cite{torsion-freepaper}. We refer the reader to the introduction in \cite{torsion-freepaper} for a more detailed explanation of this argument, and how the affine grassmannian comes in.

\begin{theorem}[{\cite[Thm. 4.8]{torsion-freepaper}}]
\label{thm: monotone}
The polynomial numerical invariant in Definition \ref{defn: numerical invariant} is strictly $\Theta$-monotone and strictly $S$-monotone on the stack $\Coh^{tf}_r(X)$.
\end{theorem}

Recall the notion of HN boundedness in Definition \ref{defn: HN boundedness}. The following theorem was proven by a concrete optimization argument in \cite{torsion-freepaper}.
\begin{theorem}[{\cite[Prop. 5.5]{torsion-freepaper}} ]
\label{thm:HNboundedness}
The polynomial numerical invariant in Definition \ref{defn: numerical invariant} on the stack $\Coh^{tf}_r(X)$ satisfies the HN boundedness condition. 
\end{theorem}
\end{subsection}

Using the $\Theta$-monotonicity condition satisfied by the stack of torsion-free coherent sheaves in Theorem \ref{thm: monotone}, plus the HN boundedness condition of the polynomial numerical invariant $\nu$ 
in Theorem \ref{thm:HNboundedness}, it follows that there exists a $\Theta$-stratification on $\Coh^{tf}_r(X)$ (Definition \ref{def:theta-filtration}), as stated by the Intrinsic GIT Theorem \ref{thm: theta stability paper theorem}(1). 

\begin{theorem}[{\cite[Thm. 5.11]{torsion-freepaper}}]
\label{thm: theta stratification torsion-free}
The polynomial numerical invariant $\nu$ defines a $\Theta$-stratification on the stack $\Coh^{tf}_r(X)$.
\end{theorem}

From the polynomial numerical invariant in Definition \ref{defn: numerical invariant} we can derive a definition of $\nu$-stability by Definition \ref{defn: semistability-nu} in terms of $\Theta$-filtrations or, by Proposition \ref{prop: filtrations of torsion free sheaves}, $\mathbb{Z}$-filtrations.  Recall the total order on polynomials in Definition \ref{defn: order_on_polynomials}.

\begin{definition}[$\nu$-semistability]\label{defn: Gieseker semistable}
Let $K \supset k$ be an algebraically closed field extension and let $\mathcal{F}$ be a torsion-free coherent sheaf on $X_K$. We say that $\mathcal{F}$ is $\nu$-semistable if for all $\mathbb{Z}$-filtrations $(\mathcal{F}_{m})_{m \in \mathbb{Z}}$ we have
\[ \sum_{m \in \mathbb{Z}} m \cdot (\overline{p}_{\mathcal{F}_{m}/\mathcal{F}_{m+1}} - \overline{p}_{\mathcal{F}}) \cdot \rk(\mathcal{F}_m/\mathcal{F}_{m+1}) \leq 0 \]
Otherwise, $\mathcal{F}$ is called $\nu$-unstable.
\end{definition}

In \cite[Prop. 5.16]{torsion-freepaper} it is shown that a torsion-free coherent sheaf $\cF$ is $\nu$-semistable  if and only if it is Gieseker semistable as in \cite{gieseker_torsion_free, maruyama_moduli, simpson-repnfundamental-I}. Therefore the stack intrinsic definition of $\nu$-stability for this polynomial  numerical invariant happens to coincide with the usual polynomial Gieseker stability for torsion-free sheaves the in literature.

Using the definition of $\Theta$-stratification, the substack of semistable torsion-free sheaves of rank $r$, denoted by $\Coh^{tf}_{r}(X)^{\nu\dash \rm{ss}}$, is an open substack of $\Coh^{tf}_{r}(X)$. 
Given a rational polynomial $P\in \mathbb{Q}[n]$, we denote by $\Coh^{tf}_{r}(X)_{P}^{\nu\dash \rm{ss}}$ the open and closed substack of $\Coh^{tf}_{r}(X)^{\nu\dash \rm{ss}}$ parametrizing sheaves whose Hilbert polynomial is $P$.

\begin{proposition}[{\cite[Prop. 5.18]{torsion-freepaper}}]
\label{prop: torsion-free bounded}
The substack $\Coh^{tf}_{r}(X)_{P}^{\nu\dash \rm{ss}}$ is quasi-compact. 
\end{proposition}

Finally, we gather Theorem \ref{thm: theta stratification torsion-free}, the $S$-monotonicity property shown in Theorem \ref{thm: monotone} and last Proposition \ref{prop: torsion-free bounded} to prove the existence of a separated good moduli space by Theorem \ref{thm: theta stability paper theorem} (2). Together with the proof of the existence part of the valuative criterion for properness in \cite[Prop. 5.20]{torsion-freepaper} we obtain an intrinsic construction of a proper good moduli space of torsion-free sheaves as in Theorem \ref{thm: theta stability paper theorem} (3). 

\begin{theorem}[{\cite[Thm. 5.21]{torsion-freepaper}}]
\label{thm: torsion-free good moduli}
The stack $\Coh^{tf}_{r}(X)_{P}^{\nu\dash \rm{ss}}$ admits a proper good moduli space.
\end{theorem}

\begin{subsection}{Leading term and Gieseker-Harder-Narasimhan filtrations}
\label{ssec: LT and GHN for torsion free}
Let $K$ be an algebraically closed field and let $\mathcal{F}$ be an unstable torsion-free coherent sheaf on $X_{K}$. 
Recall from Definition \ref{defn: hn-filtration} that a $\Theta$-filtration $f$, with associated underlying filtration $(\mathcal{F}_m)_{m \in \mathbb{Z}}$, is called a $\nu$-Harder-Narasimhan filtration if the polynomial numerical invariant $\nu(f)$ achieves its maximal value.
Having a $\Theta$-stratification on $\Coh^{tf}_r(X)$ allows us to define a unique $\nu$-HN filtration maximizing $\nu$, as seen at the end of Section \ref{ssec: theta stratifications}, for any field valued point of the stack.

\begin{proposition} \label{prop: canonical filtrations}
Let $K \supset k$ be an arbitrary field extension and let $\mathcal{F}$ be an unstable torsion-free coherent sheaf on $X_{K}$. Then, $\mathcal{F}$ admits a $\nu$-Harder-Narasimhan filtration $(\mathcal{F}_{m})_{m \in \mathbb{Z}}$ which is uniquely determined up to scaling all of the indexes by a constant rational number.  We call this filtration $(\mathcal{F}_m)_{m \in \mathbb{Z}}$ the leading term HN filtration of $\mathcal{F}$.
\end{proposition}

Depending on whether we use slope
semistability by Mumford \cite{mumford-projetive-invariants}, or Gieseker semistability \cite{gieseker_torsion_free, simpson-repnfundamental-I}, we can consider two
different filtrations, the slope Harder-Narasimhan filtration (slope HN) and the 
Gieseker Harder-Narasimhan filtration (Gieseker HN). Other intermediate polynomial stability conditions can be considered as in \cite[Thm. 1.6.7]{huybrechts.lehn}.

The leading term HN filtration of an unstable torsion-free coherent sheaf $\mathcal{F}$ from Proposition \ref{prop: canonical filtrations} corresponds to an intermediate Harder-Narasimhan filtration where the sequence of Hilbert polynomials of the quotients find their differences not at the top-degree coefficient (as in the slope HN filtration), nor at any degree coefficient (as in the Gieseker HN filtration), but at an intermediate degree coefficient, from where the filtration takes its name. Let us explain this situation. 

Start with the Gieseker HN filtration of $\mathcal{F}$ (c.f. \cite[Thm. 1.3.4]{huybrechts.lehn}):
\[ 0 = \mathcal{F}^{\text{HN}}_{0} \subset \mathcal{F}^{\text{HN}}_{1} \subset \mathcal{F}^{\text{HN}}_{2} \subset \cdots \subset \mathcal{F}^{\text{HN}}_{n} = \mathcal{F} \]
which verifies 
\[ \overline{p}_{\mathcal{F}_1^{\text{HN}}} > \overline{p}_{\mathcal{F}_2^{\text{HN}} / \mathcal{F}_1^{\text{HN}}} > \cdots > \overline{p}_{\mathcal{F}_n^{\text{HN}} / \mathcal{F}_{n-1}^{\text{HN}}} \]
There exists a minimal index $i$  such that for all $h>i$ we have that all $h^{th}$ slopes are equal
\[\widehat{\mu}_{h}(\mathcal{F}_1^{\text{HN}}) = \widehat{\mu}_{h}(\mathcal{F}_2^{\text{HN}} / \mathcal{F}_1^{\text{HN}}) = \cdots = \widehat{\mu}_{h}(\mathcal{F}_n^{\text{HN}} / \mathcal{F}_{n-1}^{\text{HN}}),\]
but $i$ is the largest index for which we find at least one strict inequality in the sequence of $i^{th}$ slopes:
\[ \widehat{\mu}_{i}(\mathcal{F}_1^{\text{HN}}) \geq \widehat{\mu}_{i}(\mathcal{F}_2^{\text{HN}} / \mathcal{F}_1^{\text{HN}}) \geq \cdots \geq \widehat{\mu}_{i}(\mathcal{F}_n^{\text{HN}} / \mathcal{F}_{n-1}^{\text{HN}}) .\]
Define indexes $h_0 < h_1 < h_2 < \cdots < h_j = n$ for the jumps of the $\widehat{\mu}_i$ sequence: let $h_0 = 0$ and
\[\{h_m\}_{m=1}^j = \left\{\; 1 \leq l \leq n  \; \left| \; \widehat{\mu}_i(\mathcal{F}_{l+1}^{\text{HN}} / \mathcal{F}_{l}^{\text{HN}}) < \widehat{\mu}_i(\mathcal{F}_{l}^{\text{HN}} / \mathcal{F}_{l-1}^{\text{HN}}) \right. \right\}\]
 and set $\mathcal{F}^{\text{l-term}}_{k} = \mathcal{F}^{\text{HN}}_{h_k}$. 
 
The filtration
\[ 0 = \mathcal{F}^{\text{l-term}}_{0} \subset \mathcal{F}^{\text{l-term}}_{1} \subset \mathcal{F}^{\text{l-term}}_{2} \subset \cdots \subset \mathcal{F}^{\text{l-term}}_{j} = \mathcal{F} \]
is called the unweighted leading term filtration of $\mathcal{F}$. We call $i$ the leading term index, which coincides with the degree of the polynomial $\nu(f)$ has degree $i$, and call $\widehat{\mu}_i(\mathcal{F})$ the leading slope. Observe that, in the case 
that $\mathcal{F}$ is slope unstable, then $i=d-1$ and hence the leading term HN filtration coincides with the slope HN filtration.

\begin{definition}[{\cite[Thm. 5.11]{torsion-freepaper} \cite[Defn. 9.2]{rho-sheaves}}]
\label{defn: canonical LT filtration}
Let $\mathcal{F}$ be a $\nu$-unstable torsion-free coherent sheaf. If the unweighted leading term HN filtration for $\mathcal{F}$ is given by
\[ 0 = \mathcal{F}^{\text{l-term}}_{0} \subset \mathcal{F}^{\text{l-term}}_{1} \subset \mathcal{F}^{\text{l-term}}_{2} \subset \cdots \subset \mathcal{F}^{\text{l-term}}_{j} = \mathcal{F} \]
with leading slope $\widehat{\mu}_{i}(\mathcal{F})$, we define the (weighted) leading term HN filtration $(\mathcal{F}^{\text{Can}})_{m \in \mathbb{Z}}$ as follows (up to scaling of the weights).

Set $\mu = \widehat{\mu}_{i}(\mathcal{F})$, $\mu_k = \widehat{\mu}_{i}\left(\mathcal{F}^{\text{l-term}}_{k} / \mathcal{F}^{\text{l-term}}_{k-1}\right)$. Choose a positive rational number $L$ such that $L\cdot \mu\in \mathbb{Z}$ and $L \cdot \mu_k \in \mathbb{Z}$ for all $k$. Then, 
\begin{enumerate}[(1)]
    \item $\mathcal{F}^{\text{Can}}_m = 0$ for $m > L (\mu_1 - \mu)$.
    \item Let $1 \leq k \leq j-1$. For each $L(\mu_k - \mu) \geq m > L (\mu_{k+1} - \mu)$, we set $\mathcal{F}^{\text{Can}}_m = \mathcal{F}^{\text{l-term}}_{k}$. 
    \item $\mathcal{F}^{\text{Can}}_m = \mathcal{F}$ for $L(\mu_{j} - \mu) \geq m$.
\end{enumerate}
This is precisely the leading term HN filtration where the $k^{th}$ jump happens at the integer $L (\mu_k - \mu)$.
\end{definition}

The weighted leading term HN filtration recovers a destabilizing HN-type filtration at the maximal degree-coefficient of the Hilbert polynomial, the leading term, yielding a sequence of decreasing leading slopes. If we take the graded object of the filtration and iterate the process we will check all possible instability terms at all lower degrees in the polynomial. This means that we will arrive at the classical Gieseker HN filtration at the end. Let us explain this procedure (c.f. \cite[Prop. 9.3]{rho-sheaves}).

Denote the (weighted) leading term HN filtration in Definition \ref{defn: canonical LT filtration} by $\left( ^{(1)}\mathcal{F}_m\right)_{m\in \mathbb{Z}}$ and the leading term index by $i_1$. Define the $1^{st}$ associated graded Levi sheaf as 
\[{\rm gr}_1(\mathcal{F}) = \left( \mathcal{F}_i^{l-term} / \mathcal{F}_{i-1}^{l-term}\right)_{i=1}^{j_1}\]
The recursion step consists of taking the leading term HN filtration of this graded object, which amounts to taking the leading term HN filtration for each $\mathcal{F}_i^{l-term} / \mathcal{F}_{i-1}^{l-term}$, with $i=1, 2, \ldots, j_1$. Note that all leading term indices are strictly smaller than $i_1$. 

This yields a lexicographic filtration indexed by $\mathbb{Z}^2$, denoted by $\left( ^{(2)}\mathcal{F}_{\overrightarrow{m}}\right)_{\overrightarrow{m}\in \mathbb{Z}^2}$, with jumps at the pairs
\[\overrightarrow{m}_1 > \overrightarrow{m}_2 > \cdots > \overrightarrow{m}_{j_2}\]
This way, the underlying unweighted filtration reads
\[ 0  \subset ^{(2)}\mathcal{F}_{\overrightarrow{m}_1}\subset ^{(2)}\mathcal{F}_{\overrightarrow{m}_2}\subset \cdots \subset ^{(2)}\mathcal{F}_{\overrightarrow{m}_{j_2}} = \mathcal{F} \]
and satisfies that the Hilbert polynomials of the graded pieces $^{(2)}\mathcal{F}_{\overrightarrow{m}_i} / ^{(2)}\mathcal{F}_{\overrightarrow{m}_{i-1}}$ are decreasing but their differences occur, at most, at a degree $i_2<i_1$, $i_2$ being the leading term index of this second leading term HN filtration.

Continuing this process (i.e. taking the leading term HN filtration of the graded object while possible), we arrive to a $q^{th}$ filtration $\left( ^{(q)}\mathcal{F}_{\overrightarrow{m}}\right)_{\overrightarrow{m}\in \mathbb{Z}^q}$ with unweighted filtration
\[ 0  \subset ^{(q)}\mathcal{F}_{\overrightarrow{m}_1}\subset ^{(q)}\mathcal{F}_{\overrightarrow{m}_2}\subset \cdots \subset ^{(q)}\mathcal{F}_{\overrightarrow{m}_{j_q}} = \mathcal{F} \]
whose graded pieces forming the $q^{th}$ associated Levi sheaf ${\rm gr}_q(\mathcal{F}) = \left( ^{(q)}\mathcal{F}_{\overrightarrow{m}_{i}}/ ^{(q)}\mathcal{F}_{\overrightarrow{m}_{i-1}}\right)$ do not admit any further filtering by any lower leading degree: they are Gieseker semistable. Therefore this final filtration is, indeed, the Gieseker HN filtration that we find in the literature \cite{simpson-repnfundamental-I, huybrechts.lehn}.

\end{subsection}
\end{section}

\begin{section}{$G$-bundles on higher dimensional varieties}

In this section we expose the theory of the ``Beyond GIT'' program for the moduli problem of $G$-bundles on higher dimensional varieties, providing a stack-intrinsic construction of the moduli space and a new notion of a Gieseker Harder-Narasimhan filtration. A more extensive treatment can be found in \cite{rho-sheaves}. 

For the rest of this section $X$ will denote a smooth projective variety of dimension $d$ over an algebraically closed field $k$ of characteristic $0$. 

\subsection{Principal $G$-bundles}

Let $G$ be an algebraic reductive group.  Let \linebreak $\rho:G\inj \GL(V)$ be a
faithful representation. A principal $G$-bundle $P$ on a scheme $X$ 
induces a vector
bundle $\cF:=P(V)$ of rank equal to $\dim V$. Using the converse construction, we
want to describe a principal bundle as a pair consisting of a vector bundle and a
reduction of structure group. 

Let $S$ be an affine $k$-scheme. Given a vector bundle $\cF$ on $X_S$ of rank
$\dim V$, we define the associated principal $\GL(V)$-bundle
\[
P_{\GL}={\underline\Iso}(\cF,V\otimes \cO_{X_S}),
\] 
with the right action of $\GL(V)$ defined as the post-composition with
the inverse of $h\in \GL(V)$:
\[
(\beta,h) \mapsto \beta\cdot h := h^{-1}\circ \beta,
\]
where we take the inverse to obtain an action on the right.

Define the scheme of reductions as the quotient
\[
{\rm Red_G(\cF)=\underline\Iso(\cF,V\otimes \cO_{X_S})/G} \longrightarrow X_S,
\]
where $G$ acts on the right via the representation $\rho$.
A reduction of structure group is given by a section of $\Red_G(\cF)$
over $X_S$.

If $\dim X>1$, then the moduli space of vector bundles does not satisfy the existence part of the valuative criterion for properness in general, because a family of vector bundles might degenerate to a
torsion-free sheaf. Likewise, if we consider principal $G$-bundles we have
to allow  for more general objects to appear in the moduli problem. These will be called principal
$\rho$-sheaves. 

Before considering the general case, let us have a look at a simpler example. A subset $U \subset X_S$ is called big if for any point $s\in X$, the fiber $U_s \subset X_s = X$
is an open set with codimension in $X$ at least two.

\begin{example}
If $G={\rm Sp}(r)$ is the symplectic group and 
$\rho$ is the standard representation, then a reduction is the same
thing as a homomorphism \linebreak $\varphi:\cF\otimes \cF \rightarrow \cO_{X_S}$ which
is bilinear, skew-symmetric and non-degenerate. 
In this case, a $\rho$-sheaf on $X$ 
will be a pair 
\[
(\cF,\varphi:\cF\otimes \cF \rightarrow \cO_X),
\]
where $\cF$ is a torsion-free sheaf and $\varphi$ is a
skew-symmetric bilinear form, which is non-degenerate when restricted to the open set
$U_{\cF}$ where $\mathcal{F}$ is locally free.
Equivalently, a $\rho$-sheaf is a principal ${\rm Sp}(r)$-bundle $P$ on a big
open set $U\subset X$, a
torsion-free sheaf $\cF$ on $X$, and an isomorphism between the vector bundle associated to $P$ 
and the restriction $\cF|_U$.
\end{example}

\begin{remark}
\label{rmk:tuples}
If we want our definition of stability to reduce to
the definition of Ramanathan when $\dim X=1$ \cite{ramanathan-stable}, then we have to require
that the representation $\rho:G \to \GL(V)$ is central (meaning that
the representation $\rho$ sends the identity component of the center
of $G$ to the center of $\GL(V)$). Otherwise, for noncentral $\rho$ our definition of semistability is stronger than Ramanathan's, see \cite[Example 10.10]{rho-sheaves} for an explanation. Notice that if $G$ is semisimple, then $\rho$ is automatically central. On the other hand, $\rho$ cannot be central if the
dimension of the center of $G$ is larger than one; this is why  
we consider representations into products of general linear groups in subsection \ref{subsect:tuples}. 
\end{remark}

\begin{subsection}{Principal $\rho$-sheaves}

%Explain that $G$-bundles on higher dimensional varieties %don't satisfy existence of the valuative criterion %properness. Introduce the moduli problem of $\rho$-sheaves.

The representation $\rho:G\inj \GL(V)$ gives a left action of $G$ on
$V$, hence a left action on $\det(V)$, and therefore $\det(V)$
acquires the structure of a right $k[G]$-comodule:
\begin{equation}
  \label{comoddet}
\det(V) \longrightarrow \det(V)\otimes k[G].
\end{equation}
Explicitly, $v \mapsto v\otimes \det(\rho(g))$.

Using the inverse $v\mapsto g^{-1}v$ we obtain a right action
of $G$ on $V$. This right action of $G$ on $V$ induces on $V$ a left $k[G]$-comodule structure $V\to
k[G]\otimes V$. We equip $V^\vee$ with the dual right
$k[G]$-comodule structure:
\begin{equation}
  \label{comodv}
c:V^\vee \longrightarrow V^\vee \otimes k[G].
\end{equation}
Explicitly, if $\{e_i\}$ is a basis of $V$ and $\{e^i\}$ is the dual
basis of $V^\vee$, then
$$
c(e^i) = \sum_j e^j \otimes (\rho(g)^{-1})^i{}_j
$$
(note that  $A\mapsto (A^{-1})^i{}_j$ is an element of $k[\GL(V)]$). 
The value that it takes on an element $A\in \GL(V)$ is the  $(i,j)$-entry of the inverse 
of $A$, expressed in the basis $\{e_i\}$.

These right $k[G]$-comodule structures \eqref{comoddet} \eqref{comodv} induce a right
$k[G]$-comodule structure on
\begin{equation}
  \label{evf}
E_V(\cF) := \Sym^\bullet(\cF\otimes_k V^\vee)\otimes_{\cO_{X_S}} 
\Sym^\bullet(\det(\cF)^\vee \otimes_k \det(V)).
\end{equation}
Therefore, the scheme $H(\cF,V):= \underline{\Spec}_{X_S}(E_V(\cF))$
(which is affine over $X_S$ and of finite presentation) has
a right action by $G$. Consider the GIT quotient
\[
H(\cF,V) {/\!\!/} G = \underline{\Spec}_{X_S}(E_V(\cF)^G),
\]
where $E_V(\cF)^G$ is the subalgebra of $G$-invariants.

The restriction $H:=H(\cF,V)|_{U}$ to the open set $U\subset X_S$ where
$\cF$ is locally free is 
\[
H= \underline\Hom(\cF|_U,V\otimes_k \cO_U) \times
\underline\Hom(\det(\cF|_{U})^\vee,\det(V^\vee)\otimes_k \cO_U),
\]
parametrizing pairs $(\alpha:\cF|_U\to V\otimes_k \cO_U,
\beta:\det(\cF|_{U})^\vee\to \det(V^\vee)\otimes_k \cO_U$. The equation
$\det(\alpha)\beta=1$ cuts out a closed subscheme which is isomorphic to
\[
\underline\Iso(\cF|_U,V\otimes_k \cO_U).
\]
Furthermore, the equation $\det(\alpha) \beta=1$ is $G$-invariant, so it
also cuts out a closed subscheme  
$(H{/\!\!/} G)_0 \subset H{/\!\!/} G$ which is isomorphic
to 
\[
\underline\Iso(\cF|_U,V\otimes_k \cO_U) /G = \Red_G(\cF|_U).
\]
\begin{definition}[{\cite[Defn. 2.14]{rho-sheaves}}]
Let $\cF$ be a torsion-free sheaf on $X$. The scheme of reductions
$\Red_G(\cF)$ is defined to be the scheme theoretic closure of 
$\Red_G(\cF|_U)$ in $H(\cF,V) {/\!\!/} G$. 
\end{definition}

The scheme $\Red_G(\cF) \to X_S$ is affine of finite type and allows to give a definition of principal $\rho$-sheaf.

\begin{definition}[{\cite[Defn. 3.1]{rho-sheaves}}]
A principal $\rho$-sheaf on $X_S$ is a pair $(\cF,\sigma)$ where
\begin{enumerate}
\item $\cF$ is a torsion-free sheaf on $X_S$.
\item $\sigma$ is a section of $\Red_G(\cF) \to X_S$.
\end{enumerate}
\end{definition}

We have the following Hartogs-like result (see \cite[Lem. 2.2]{rho-sheaves}).

\begin{proposition}[{\cite[Prop. 2.18]{rho-sheaves}}]
\label{hartoglike}
Let $U\subset X_S$ be any big open set where $\cF$ is locally free (note that such $U$ always exists \cite[Lemma 2.5]{rho-sheaves}).
There is a bijection between the sections of the morphism 
$\Red_G(\cF) \to X_S$ and the sections on the restriction 
$\Red_G(\cF|_U) \to U$.
\end{proposition}

This allows us to give an equivalent definition of $\rho$-sheaf.

\begin{lemma}
Giving a principal $\rho$-sheaf on $X_S$ is the same thing as giving a triple $(P,\cF,\psi)$ where
\begin{enumerate}
\item $\cF$ is a torsion-free sheaf on $X_S$, and $U$ is a big open subset where $\cF$ is locally free (we may take $U$ to be the maximal such open).
\item $P$ is a principal $G$-bundle on $U$.
\item $\psi$ is an isomorphism between the vector bundle $P(\rho)$ 
associated to $\rho$, and the restriction $\cF|_U$
\[
\psi:P(\rho)\cong \cF|_U .
\]
\end{enumerate}
\end{lemma}

Denote by $\Bun_{\rho}(X)$ the stack of principal $\rho$-sheaves on $X$.

\begin{definition} The algebraic stack of principal $\rho$-sheaves on $X$ is the pseudofunctor defined by
\[
T \; \mapsto \; \Bun_{\rho}(X)(T) \; = \;   \big[\; \text{groupoid of $\rho$-sheaves $(\cF,\sigma)$
  on $X_T$} \;\big] .
\]
\end{definition}

There is a forgetful morphism $\Bun_{\rho}(X)\rightarrow \Coh_{r}^{tf}(X)$
which is schematic, affine and of finite type (\cite[Prop 3.6]{rho-sheaves}). This allows to relate the constructions for $\rho$-sheaves to those of section \ref{sec: torsion free} for torsion-free coherent sheaves. 

\end{subsection}

\begin{subsection}{Moduli space of $\rho$-sheaves}
\label{subs:modulirho}
%Define Gieseker stability. State the theorem about existence of moduli %space. (Possibly: Explain some interesting examples of stability, and %checking on maximal parabolics).

Here we introduce a polynomial numerical invariant on the stack $\Bun_{\rho}(X)$ of principal $\rho$-sheaves, which yields the notion of Gieseker stability. Once we prove that the requirements of Theorem \ref{thm: theta stability paper theorem} are met, it will be shown that 
the substack of semistable $\rho$-sheaves admits a good moduli space.

Let $\lambda:\mathbb{G}_m\to G$ be a  one-parameter subgroup . We define
a parabolic subgroup $P_\lambda\subset G$ associated to $\lambda$ by
\[
P_{\lambda} = \left\{ \, g \in G \, \mid \, \lim \limits_{x \to 0} \, \lambda(x) \,g\, \lambda(x)^{-1} \;\text{exists} \,\right\}
\]
whose Levi subgroup is
\[
L_{\lambda} = \left\{ \, g \in G \, \mid \, \lim \limits_{x
       \to 0} \, \lambda(x) \,g\, \lambda(x)^{-1} = g \,\right\}.
\]
Composing with the representation $\rho$ we obtain a parabolic subgroup
of $\GL(V)$:
$$
P^{\GL(V)}_{\lambda} = \left\{ \, g \in \GL(V) \, \mid \, \lim \limits_{x \to 0} \, \rho\circ\lambda(x) \,g\, (\rho\circ\lambda(x))^{-1} \;\text{exists} \,\right\}
$$
and it is clear that $P_\lambda=P^{\GL(V)}_\lambda\cap G$.

Let $(\cF,\sigma)$ be a principal $\rho$-sheaf on $X$. Let $P\to U\subset X$
be the associated principal $G$-bundle on the open set $U$ where $\cF$
is locally free, and let 
$\mathcal{G}_\lambda\to U'\subset U$ be a reduction of structure
group of $P$ to $P_{\lambda}$ on a big open set $U'\subset X$. The inclusion
$P_{\lambda}\subset P^{\GL(V)}_{\lambda}$ gives a principal
$P^{\GL(V)}_\lambda$-bundle $P'$ on $U'$, which is a reduction of structure group
to the parabolic $P^{\GL(V)}_\lambda$ of the principal $\GL(V)$-bundle
associated to the vector bundle $\cF|_{U'}$. 

Giving a reduction to a parabolic of the linear group is the same thing as giving a filtration of the vector bundle $\cF|_{U'}$, and such a filtration extends uniquely over $X_S$ to a
filtration of $\cF$ by torsion-free sheaves (with the condition 
that successive quotients of the filtration are also torsion-free). 
Moreover, we use the weights of
the one-parameter subgroup obtained by composition $\rho\circ\lambda:
\mathbb{G}_m\to \GL(V)$ to obtain a $\mathbb{Z}$-weighted filtration
$(\cF_m)_{m\in \mathbb{Z}}$ of $\cF$
as in Proposition \ref{prop: filtrations of torsion free sheaves}.

Recall from Definition \ref{def:theta-filtration} that a $\Theta$-filtration of a point $p$ in a stack $\cM$ 
is a morphism $f:\Theta=[\mathbb{A}^1/\mathbb{G}_m]\to \cM$ together with
an isomorphism  $f(1)\simeq p$. Any \linebreak $\Theta$-filtrations $f:\Theta\rightarrow \Bun_{\rho}(X)$ of a $\rho$-sheaf $(\mathcal{F}, \sigma)$ can be composed with the forgetful morphism $\Bun_{\rho}(X)\rightarrow \Coh_{r}^{tf}(X)$ to obtain a $\mathbb{Z}$-weighted filtration of $\mathcal{F}$ satisfying a series of compatibility conditions. This way, a $\Theta$-filtration of a $\rho$-sheaf $(\mathcal{F}, \sigma)$ can be understood as a $\mathbb{Z}$-weighted filtration
$(\cF_m)_{m\in \mathbb{Z}}$ of $\cF$
(as in Proposition \ref{prop: filtrations of torsion free sheaves}) which is associated to a reduction of structure group (\cite[Prop. 3.12 and 3.16]{rho-sheaves}). See {\cite[Sections 3.2 and 3.3]{rho-sheaves}} for details of this construction.

In a similar way to section \ref{ssec: tf sheaves}, we can define a family of line bundles \cite[Defn. 4.1]{rho-sheaves} and a rational quadratic norm \cite[Defn. 4.3]{rho-sheaves} to define a polynomial numerical invariant (\cite[Prop. 4.4]{rho-sheaves}):
\[
\nu(f)= \sqrt{A_d} 
\frac{\sum_{m \in \mathbb{Z}} m \cdot
  (\overline{p}_{\cF_{m}/\cF_{m+1}} - \overline{p}_{\cF}) \cdot
  \rk(\cF_m/\cF_{m+1})}
{\sqrt{\sum_{m\in \mathbb{Z}} m^2 \rk(\cF_m/\cF_{m+1})}}.
\]
Note that this is the same expression as in Definition \ref{defn: numerical invariant} for torsion-free sheaves, but for a $\rho$-sheaf $(\mathcal{F}, \sigma)$ it accounts just for $\Theta$-filtrations of $\Bun_{\rho}(X)$.

The polynomial numerical invariant $\nu$ defined on the stack $\Bun_{\rho}(X)$ can be shown to be strictly $\Theta$-monotone and strictly $S$-monotone (\cite[Prop. 4.16]{rho-sheaves}) and to satisfy the HN boundedness condition (\cite[Prop. 4.26]{rho-sheaves}). Then it defines a $\Theta$-stratification on $\Bun_{\rho}(X)$ by Theorem \ref{thm: theta stability paper theorem}(1) (\cite[Thm. 4.27]{rho-sheaves}). This allows to define a notion of $\nu$-stability which is an open condition.

\begin{definition}[{\cite[Defn. 4.28]{rho-sheaves}}]
\label{defn:ssrhosheaves}
A principal $\rho$-sheaf $(\cF,\sigma)$ is $\nu$-semistable if for all
one-parameter subgroups $\lambda:\mathbb{G}_m\to G$ and all reductions
as above, we have
\[ \sum_{m \in \mathbb{Z}} m \cdot (\overline{p}_{\cF_{m}/\cF_{m+1}} - \overline{p}_{\cF}) \cdot \rk(\cF_m/\cF_{m+1}) \leq 0 .\]
\end{definition}

This definition coincides with the semistability notion in \cite{schmitt.singular, gomezsols.principalsheaves}.
The substack $\Bun_{\rho}(X)^{\nu\dash \rm{ss}}$ of semistable $\rho$-sheaves satisfies the existence part of the valuative criterion for properness in \cite[Thm. 5.8]{rho-sheaves}. After fixing the degree $\vartheta$ of a $\rho$-sheaf $(\mathcal{F}, \sigma)$ (\cite[Defn. 6.3]{rho-sheaves}) and the Hilbert polynomial $P$ of $\mathcal{F}$, the substacks $\Bun_{\rho}(X)^{\nu\dash \rm{ss},P}_{\vartheta}$ are quasi-compact (\cite[Prop. 6.6]{rho-sheaves}, then by Theorem \ref{thm: theta stability paper theorem} (2,3) they do admit a proper good moduli space (\cite[Thm. 6.8]{rho-sheaves}). This is a stack intrinsic alternative to the constructions in \cite{schmitt.singular, gomezsols.principalsheaves}.

\end{subsection}

\begin{subsection}{The scheme of $G$-reductions for products of representations}
\label{subsect:tuples}

As we explained in Remark \ref{rmk:tuples}, if we want to consider 
reductive groups that are not semisimple, we need to deal with
tuples of representations. In section \ref{ssec: GHN for rho sheaves} we will observe that when constructing the Gieseker HN filtration, even if we start with a semisimple
group, the Levi factors we need to go through make necessary to consider non-semisimple groups. 

Let $V^\bullet=(V^i){}_{i=1}^b$ be a tuple of finite dimensional vector spaces. 
We consider a representation into the tuple
\[
\rho=(\rho_1,\ldots,\rho_b): 
G \longrightarrow \GL(V^1)\times \cdots \times \GL(V^b).
\]
Let $\cF^\bullet=(\cF^i){}_{i=1}^b$ be a tuple of torsion-free sheaves on $X_S$.
Using the construction \eqref{evf}, we consider the product scheme
\[
H(\cF^\bullet,V^\bullet):= \prod_{i=1}^b H(\cF^i,V^i).
\]
The group $G$ has a right action on this product, induced by the right
action on each factor defined by $\rho_i:G\to \GL(V^i)$. Let
$H(\cF^\bullet,V^\bullet) {/\!\!/} G$ be the GIT quotient.  
Let $U$ be an open set where all sheaves $\cF^i$ are locally free. 
On each factor $H(\cF^i|_U,V^i)$ we have an equation
$\alpha_i\det(\beta_i)=1$ as before, cutting out a divisor isomorphic
to the scheme of isomorphisms $\Iso(\cF^i|_U,V^i\otimes \cO_U)$. Finally, we define $\Red_G(\cF^\bullet)$
as the scheme theoretic closure of the subscheme
\[
\prod \Iso(\cF^i|_U,V^i\otimes \cO_U) \cong \prod (H(\cF^i|_U,V^i))_{\alpha^i\det(\beta^i)=1}\inj H(\cF^\bullet,V^\bullet) {/\!\!/} G.
\]

\begin{definition}[{\cite[Defn. 7.1]{rho-sheaves}}]
A principal $\rho$-sheaf on $X_S$ is a pair $(\cF^\bullet,\sigma)$ where

\begin{enumerate}
\item $\cF^\bullet$ is a tuple of torsion-free sheaves on $X_S$.
\item $\sigma$ is a section of $\Red_G(\cF^\bullet) \to X_S$.
\end{enumerate}
\end{definition}

As in the case of a single representation, we denote by $\Bun_\rho(X)$ the stack of $\rho$-sheaves. 

For any big open set $U$ where all sheaves $\cF^i$ are locally free,
the section $\sigma$ gives a reduction of structure group, from the
principal $\GL(V^1)\times \cdots \times \GL(V^b)$-bundle associated to the
sum of vector bundles $\cF^1|_U\oplus \cdots \oplus \cF^b|_U$, to a
principal $G$-bundle $P\to U$ on $U$. Therefore, using the Hartogs-like
Proposition \ref{hartoglike}, we have the following
\begin{lemma}
Giving a principal $\rho$-sheaf is the same thing as giving a triple
$(\cF^\bullet, P, \psi^\bullet)$ where
\begin{enumerate}
\item $P$ is a principal $G$-bundle on a big open set $U$.
\item $\cF^\bullet=(\cF^i)_{i=1}^b$ is a tuple of $b$ torsion-free
  sheaves, all of them locally free on $U$.
\item $\psi^\bullet=(\psi^i)_{i=1}^b$ is a collection of isomorphisms
$$
\psi^i:P(\rho_i)\cong \cF^i|_U\; .
$$
\end{enumerate}
\end{lemma}

Finally, we can extend to this case the notion of stability. 

\begin{definition}
\label{defn:ssrhosheavesproduct}
A principal $\rho$-sheaf $(\cF^\bullet,\sigma)$ on $X$ is $\nu$-semistable if for all
one-parameter subgroups $\lambda: \mathbb{G}_m\to G$ and all reductions
as described in subsection \ref{subs:modulirho}, we have
\[ \sum_{i=1}^{b} \sum_{m \in \mathbb{Z}} m \cdot (\overline{p}_{\cF^i_{m}/\cF^i_{m+1}} - \overline{p}_{\cF^i}) \cdot \rk(\cF^i_m/\cF^i_{m+1}) \leq 0 \]
\end{definition}

\end{subsection}

\begin{subsection}{Gieseker-Harder-Narasimhan filtration for $\rho$-sheaves}
\label{ssec: GHN for rho sheaves}

A unstable torsion-free sheaf has an associated canonical filtration,
called the Harder-Narasimhan filtration, showing that any torsion-free
sheaf can be constructed uniquely as a successive extension of
semistable torsion-free sheaves \cite[Thm. 1.3.4]{huybrechts.lehn}. As we mentioned in subsection \ref{ssec: LT and GHN for torsion free}, this can be the slope HN or the Gieseker HN filtration, depending on the stability condition we are using. 

In the case of principal bundles, there is already  in \cite{behrend-thesis} a notion of
Harder-Narasimhan reduction if we use Ramanathan's definition of
semistability (which is the direct generalization of Mumford's slope condition in \cite{mumford-projetive-invariants}). This is generalized to higher dimensional varieties in \cite{anchouche-hassan-biswas}. 
However, if we use Gieseker's notion \cite{gieseker_torsion_free} of semistability considering all terms of the Hilbert polynomial, the
right notion of Harder-Narasimhan filtration had not been defined
until the arrival of the ``beyond GIT'' theory \cite{halpernleistner2021structure}. The notion of semistability for $\rho$-sheaves (Definitions
\ref{defn:ssrhosheaves} and \ref{defn:ssrhosheavesproduct}) can be
obtained through the techniques of the ``beyond GIT'' program, 
in terms of $\Theta$-filtrations directly from the stack picture. The advantage of this approach
is that it provides a way to define a correct Gieseker Harder-Narasimhan filtration for moduli problems with a polynomial stability condition.

Let $(\mathcal{F}, \sigma)$ be  a $\rho$-sheaf and let $f:\Theta\rightarrow \Bun_{\rho}(X)$ be a  $\Theta$-filtration of $(\mathcal{F}, \sigma)$. Recall from subsection \ref{subs:modulirho} the definition of a polynomial numerical invariant which,  when considering representations into tuples $V^\bullet$ as in subsection \ref{subsect:tuples}, is written as (\cite[Section 7]{rho-sheaves}):
\[
\nu(f)= \sqrt{A_d} 
\frac{\sum_{i=1}^b\sum_{m \in \mathbb{Z}} m \cdot
  (\overline{p}_{\cF^i_{m}/\cF^i_{m+1}} - \overline{p}_{\cF^i}) \cdot
  \rk(\cF^i_m/\cF^i_{m+1})}
{\sqrt{\sum_{i=1}^b\sum_{m\in \mathbb{Z}} m^2 \rk(\cF^i_m/\cF^i_{m+1})}}.
\]
We can rephrase (Definition \ref{def:theta-stability}) the condition of $\nu$-semistability by saying that a
$\rho$-sheaf $(\cF,\sigma)$ is $\nu$-semistable if and only if for all $\Theta$-filtrations
$f:\Theta\to\Bun_\rho(X)$ with $f(1)\simeq (\cF,\sigma)$ we have
$\nu(f)\leq 0$.

Let $(\cF,\sigma)$ be a $\nu$-unstable $\rho$-sheaf. 
Given that $\Bun_{\rho}(X)$ is endowed with a \linebreak $\Theta$-stratification it is shown in \cite[Prop. 4.31]{rho-sheaves} that there is a unique\footnote{It is unique up to multiplication of all the indexes $m\in \mathbb{Z}$ by a scalar. This corresponds to reparametrization of the one-parameter subgroup by $z\mapsto z^a$} \linebreak $\Theta$-filtration $f$ giving
the maximum value of $\nu(f)$. This is a $\nu$-HN filtration of the point $(\mathcal{F}, \sigma)$ in the stack $\Bun_{\rho}(X)$ as in Definition \ref{defn: hn-filtration}.

Associated to this maximal $\Theta$-filtration $f$, we have 
a reduction of structure group to a parabolic subgroup
$P_{\lambda_1}\subset G$ for some $\lambda_1:\mathbb{G}_m\to G$, and
a $\mathbb{Z}$-filtration $\left( \cF_m\right)_{m\in \mathbb{Z}}$ of $\cF$ satisfying the compatibility conditions described in \cite[Section 3.2]{rho-sheaves}. By the same reasons explained in subsection \ref{ssec: LT and GHN for torsion free} this is called the leading term HN filtration, which is the version for $\rho$-sheaves of the filtration described there for torsion-free coherent sheaves.
In order to obtain a filtration of Gieseker type, meaning a filtration being able to compare all terms of Hilbert polynomials and whose graded object is semistable, we need to reproduce the algorithm described in that subsection \ref{ssec: LT and GHN for torsion free} translated to the stack $\Bun_{\rho}(X)$. 

Using the projection $P_{\lambda_1}\to L_{\lambda_1}$ from the parabolic group to its Levi factor, this reduction induces a principal $\rho_{\lambda_1}$-sheaf with structure group
$L_{\lambda_1}$, the Levi factor of $P_{\lambda_1}$.
The representation of this Levi is of the form 
\[\rho_{\lambda_1}:L_{\lambda_1}\to \GL(V^1)\times\cdots\times \GL(V^b),\]
and the tuple of sheaves $\cF'{}^\bullet$ is given by the successive
quotients of the filtration $\left( \cF_m\right)_{m\in \mathbb{Z}}$, i.e., $\cF'{}^i=\cF_i/\cF_{i+1}$. However, this Levi might not be
semistable yet, therefore we need to further take the leading term HN filtration of this $\rho_{\lambda_1}$-sheaf. In that case, we obtain a new reduction
to a parabolic $P_{\lambda_2}\subset L_{\lambda_1}\subset G$ and repeat this process until we obtain a reduction whose associated Levi sheaf
is semistable. The final associated filtration will be called the 
Gieseker HN filtration (see \cite[Section 8.3]{rho-sheaves}. 

\end{subsection}

\begin{subsection}{Example of a Gieseker Harder-Narasimhan filtration}

We illustrate the process of producing the Gieseker HN filtration with a concrete numerical example.

Let $G={\rm Sp}(6)$ with $\rho:{\rm Sp}(6)\to \GL(6)$ the standard
representation. Let $X=\mathbb{P}^3$. Let $p$ and $C$ be,
respectively, a given fixed point and a fixed line in $\mathbb{P}^3$. Let
$\mathcal{I}_p$ and $\mathcal{I}_C$ be their respective ideal
sheaves. Set
$$
\cF=\mathcal{I}_C \oplus (\mathcal{I}_p \oplus \mathcal{I}_p ) 
\oplus (\cO_X\oplus \cO_X) \oplus \cO_X
$$
and let $\varphi:\cF\otimes \cF\to \cO_X$ be the symplectic form
defined by the matrix
$$
\varphi=
\left(
  \begin{array}{ccccccc}
  0 & \! \cdots\! &    &  &  & 1 \\
  \vdots &     &  &  & 1 &  \\
   &     &  & 1 &  &  \\
   &   & -1 &  &    &  \\
   &  -1 &  &  &    &  \vdots \\
 - 1 &   &  &  &    \! \cdots\!  &  0\\
  \end{array}
\right)
$$
The Hilbert polynomials are
  \begin{align*}
P_{\mathcal{O}}(n)&=\frac{1}{6}(n^3+6n^2+10n+3)\, ,\\
P_{\mathcal{I}_p}(n)&=\frac{1}{6}(n^3+6n^2+10n-3)\, ,\\
P_{\mathcal{I}_C}(n)&=\frac{1}{6}(n^3+6n^2  -\; 6n-3)\, .
  \end{align*}
We will consider one-parameter subgroups of ${\rm Sp}(6)$ of the form
$$
f(t)=\operatorname{diag}\operatorname(t^a,t^b,t^c,t^{-c},t^{-b},t^{-a})
$$
with $a\geq b\geq c\geq 0$. The associated filtration $\cF_m$
satisfies
$\mathcal{F}_m^\perp=\mathcal{F}_{-m+1}^{}$, where the orthogonal of a
subsheaf $\mathcal{F}_i\subset \mathcal{F}$ is given by
$$
\mathcal{F}_i{}^\perp :=\ker  \big(\mathcal{F}\stackrel{\varphi}{\longrightarrow}\mathcal{F}^\vee \longrightarrow \mathcal{F}_i{}^\vee\big ).
$$
The bilinear form gives a reduction of structure group of
$\cF|_{\mathbb{P}^3\setminus {(C\cup p)}}$ to
${\rm Sp}(6)$, so the data $(\cF,\varphi)$ defines a principal $\rho$-sheaf.
It can be checked that it is $\nu$-unstable (Definition \ref{defn:ssrhosheaves}), and the one-parameter subgroup associated to the leading term HN filtration maximizing $\nu$
is $\lambda_1=\operatorname{diag}(t,1,1,1,1,t^{-1})$. The
filtration of sheaves $\left(^{(1)}\mathcal{F}_m\right)_{m\in\mathbb{Z}}$ associated to $\lambda_1$ is
$$
\overbrace{\mathcal{F}}^{\mathcal{F}_{-1}}  \quad\supset\quad
\overbrace
{\mathcal{I}_p \oplus \mathcal{I}_p  \oplus 
 \mathcal{O} \oplus \mathcal{O}}
^{\mathcal{F}_0}
\quad\supset\quad
\overbrace{\mathcal{O}}^{\mathcal{F}_{1}}
$$
We have the associated Levi $L_{\lambda_1}=\GL(1)\times {\rm Sp}(4)$ and
the representation $\rho_{\lambda_1}$, which is given by 
$$
\begin{array}{rccc}
\rho_{\lambda_1}:&\GL(1)\times {\rm Sp}(4)& \longrightarrow &\GL(1)\times \GL(4) \times \GL(1)\\
& (t,g) & \mapsto & (t,g,t^{-1})\\
\end{array}
.
$$
The $1^{st}$ associated Levi graded
sheaf ${\rm gr}_1(\mathcal{F})$ is
$$
\overbrace{\mathcal{I}_L}^{\mathcal{F}^{-1}:= \mathcal{F}_{-1}/\mathcal{F}_{0}}
\quad\oplus\quad 
\overbrace
{\mathcal{I}_p \oplus \mathcal{I}_p \oplus \mathcal{O} \oplus \mathcal{O} }
^{\mathcal{F}^0 := \mathcal{F}_{0}/\mathcal{F}_{-1}}
\quad\oplus\quad 
\overbrace{\mathcal{O}}^{\mathcal{F}^{1} := \mathcal{F}_{1}}.
$$
The principal $\rho_1$-sheaf defined by this is not $\nu$-semistable, then we take its leading term HN filtration having associated one-parameter
subgroup 
$$
\begin{array}{rccc}
\lambda_2:&\mathbb{G}_m& \longrightarrow &\GL(1)\times {\rm Sp(4)}\\
& t & \mapsto &  \big( 1  ,  (t,t,t^{-1},t^{-1})\big)\\
\end{array}
$$
The corresponding $2^{nd}$ filtration $\left(^{(2)}\mathcal{F}_{(m_1, m_2)}\right)_{(m_1,m_2)\in\mathbb{Z}^2}$ is
\begin{equation}
  \label{eq:ghnsymplectic}
\cF 
 \quad\supset \quad
 (\mathcal{I}_p \oplus \mathcal{I}_p ) 
\oplus (\cO_X\oplus \cO_X) \oplus \cO_X 
\quad\supset\quad
 (\cO_X\oplus \cO_X) \oplus \cO_X 
\quad\supset\quad
\cO_X 
\end{equation}
whose $2^{nd}$ associated Levi graded sheaf ${\rm gr}_2(\mathcal{F})$ is
$$
\mathcal{I}_C \quad\oplus\quad (\mathcal{I}_p \oplus \mathcal{I}_p ) 
\quad\oplus\quad (\cO_X\oplus \cO_X) 
\quad\oplus\quad \cO_X
$$
The second Levi factor $L_{\lambda_2}$ is $\GL(1)\times \GL(2)$
(note that the embedding $\GL(2)\subset {\rm Sp}(4)$ is given by 
$A\mapsto (A,A^{-1})$) and the principal $\rho_2$-sheaf is $\nu$-semistable.
Since this is semistable, the process stops, and the
Gieseker Harder-Narasimhan reduction has \eqref{eq:ghnsymplectic} as
associated filtration.

\end{subsection}
\end{section}
\medskip

\medskip

\noindent \textbf{Acknowledgements:} We thank Andr\'es Ib\'a\~{n}ez N\'u\~{n}ez, Sveta Makarova and Nitin Nitsure, as well as the anonymous referee, for helpful comments on the manuscript. We also thank the VBAC community lead by Peter Newstead and the editors of this contributed volume for their continuous work. 

This work is supported
by grants CEX2019-000904-S and PID2019-108936GB-C21 (funded by MCIN/AEI/ 10.13039/501100011033).

\bibliographystyle{amsalpha}

\bibliography{survey-ams-beyond-git.bib}

\end{document}